\documentclass{article}

\usepackage{amsmath,amsfonts,amssymb}
\usepackage{amsthm}
\usepackage{ifthen}
\usepackage{makeidx}
\usepackage{url,hyperref}
\usepackage{multicol}
\usepackage{graphicx}
\usepackage{multirow}

\title{Conjugacy classes and rational period functions for the Hecke groups}

\author{Wendell Ressler}

\date{\today}

	\newcommand{\Z}{\mathbb{Z}}			
	\newcommand{\R}{\mathbb{R}}			
   
\theoremstyle{plain}
	
	\newtheorem{corollary}{Corollary}

	\newtheorem{lemma}{Lemma}

    \newtheorem{theorem}{Theorem}

\theoremstyle{definition}

	\newtheorem{example}{Example}
	   
\theoremstyle{remark}

\begin{document}

\maketitle

\begin{abstract}
We establish a one-to-one correspondence
between conjugacy classes of any Hecke group
and
irreducible systems of poles
of rational period functions
for automorphic integrals 
on the same group.
We use this correspondence to
construct irreducible systems of poles
and to count poles.
We characterize Hecke-conjugation 
and Hecke-symmetry
for poles of rational period functions
in terms of the transpose 
of matrices in conjugacy classes.
We construct new rational period functions
and 
families of rational period functions.

\vspace{\baselineskip}    

\vspace{\baselineskip}    

\end{abstract}


\section{Introduction}

Marvin Knopp first 
defined and studied rational period functions (RPFs)
for automorphic integrals 
\cite{MR0344454}
and gave the first example
of an RPF
with nonzero poles 
\cite{MR0485700}.

Knopp \cite{MR623004}, 
Hawkins \cite{Hawkins_manuscript}, 
and Choie and Parson 
\cite{MR1045397,MR1103673} 
took the main steps toward an explicit characterization of RPFs on 
the modular group \(\Gamma(1)\).
An important tool for this work
was Hawkins' idea of studying
irreducible systems of poles (ISPs)
for RPFs.
Ash \cite{MR980298} 
gave an abstract
characterization, 
and then
Choie and Zagier \cite{MR1210514} and
Parson \cite{MR1210515} provided 
a more explicit characterization 
of the RPFs on \(\Gamma(1)\).
The explicit characterizations use continued fractions 
to establish a connection between the poles of 
RPFs and binary quadratic forms.

Schmidt \cite{MR1378572} generalized Ash's work
and gave an abstract characterization of RPFs 
on any finitely generated Fuchsian group of the 
first kind with parabolic elements, 
a class of groups which includes the Hecke groups.
Schmidt \cite{MR1219337} 
and Schmidt and Sheingorn \cite{MR1362251} 
gave explicit descriptions of certain RPFs 
on the Hecke groups 
using generalizations of the 
classical continued fractions and binary quadratic forms.
This author continued that work 
in \cite{MR1876701}, \cite{MR2582982}, and \cite{MR3574636}
and has explicitly characterized 
RPFs on the Hecke groups for two cases.
The question of a complete explicit
characterization remains open.

In this paper we
construct new RPFs 
on Hecke groups,
including ones
in every case
for the first four Hecke groups.
We also construct families of RPFs
across all 
(or for some families all but several) Hecke groups.
Our primary tool is
an explicit correspondence between 
poles of an RPF
and conjugacy classes of its Hecke group.
We show that by
using this correspondence 
we may construct all ISPs
of a given size
from products of generators for corresponding
conjugacy classes.
We give a formula
for the number of such products
and thus for the number of corresponding poles.
We also 
use the correspondence 
to conjugacy classes
to characterize
Hecke-conjugation and Hecke-symmetry
for poles of rational period functions.


\section{Background}

We first summarize background ideas 
and definitions necessary
to work with 
conjugacy classes and rational period functions.
More details can be found in
\cite{MR1876701},
\cite{MR2582982}, 
\cite{MR3078226}, and
\cite{MR3574636}.


\subsection{Hecke groups and Hecke symmetry}

Let \(S = S_{\lambda} 
   = \bigl( \begin{smallmatrix} 1&\lambda\\ 0&1 \end{smallmatrix} \bigr) \),
\(T = \bigl( \begin{smallmatrix} 0&-1\\ 1&0 \end{smallmatrix} \bigr)\),
and \(I = \bigl( \begin{smallmatrix} 1&0\\ 0&1 \end{smallmatrix} \bigr)\),
where \( \lambda \) is a positive real number.
Put \(G(\lambda)=\langle S,T \rangle / \{\pm I\} \subseteq \mathrm{PSL}(2,\R) \).
Erich Hecke \cite{MR1513069} showed that the  
values of \(\lambda\) for which \(G(\lambda)\) is discrete 
are 
\[\lambda = \lambda_{p} = 2\cos(\pi/p), \] 
for \(p=3, 4, 5, \dotsc\),
and for \( \lambda \geq 2 \).
The groups with \( \lambda = \lambda_{p} \), for \( p \geq 3 \)
have come to be known as the \emph{Hecke 
groups};
we denote them by
\(G_{p}=G(\lambda_{p})\) for \(p\geq3\).
The first several of these Hecke groups are 
\(G_{3}=G(1)=\Gamma(1)\) (the modular group), 
\( G_{4}=G(\sqrt{2}) \), 
\(G_{5}=G\left( \frac{1+\sqrt{5}}{2} \right)\), and
\( G_{6}=G(\sqrt{3}) \).

A Hecke group 
\( G_{p} \)
is the free product of 
the cyclic group of order \( 2 \)
generated by \( T \)
and
the cyclic group of order \( p \)
generated by 
\( U = U_{\lambda} = S_{\lambda}T 
= \bigl( \begin{smallmatrix} \lambda & -1\\ 1 & 0 \end{smallmatrix} \bigr) \),
so the group relations of \( G_{p} \) are
\( T^{2}=U^{p}=I\).

Elements of \( G_{p} \)
act on the Riemann sphere 
as M\"{o}bius transformations.
An element 
\(M 
= \bigl( \begin{smallmatrix} a&b\\ c&d \end{smallmatrix} \bigr) 
\in G_{p}\) 
is 
\emph{hyperbolic} if \( \vert a+d \vert >2 \),
\emph{parabolic} if \( \vert a+d \vert =2 \),
and \emph{elliptic} if \( \vert a+d \vert <2 \). 
We designate fixed points accordingly.
Hyperbolic M\"{o}bius transformations each have 
two distinct real fixed points.

The action of \( G_{p} \)
induces an equivalence relation 
for points 
on the Riemann sphere.
We say that 
\( z_{1} \) and \( z_{2} \) are 
\( G_{p} \)-\emph{equivalent}
if there exists 
\( M \in G_{p} \) such that
\( z_{2} = Mz_{1} \).
Equivalence classes contain either all fixed 
points of the same kind, or no fixed points.

The \emph{stabilizer} in \( G_{p} \) of a complex number \( z \),
\( \mathrm{stab}(z) = \{M \in G_{p} \mid Mz=z\} \),
is a cyclic subgroup of \( G_{p} \).
We define the \emph{Hecke-conjugate} 
of any hyperbolic fixed point \( \alpha \)
to be the other fixed point of the elements in its stabilizer
and denote it by \( \alpha^{\prime} \).
If \( R \) is a set of hyperbolic fixed points of 
\( G_{p} \) we write
\( R^{\prime} = \{x^{\prime} \mid x \in R\} \).
We say that a set \( R \) has 
\emph{Hecke-symmetry}
if \( R = R^{\prime} \).
We observe that 
\( R \cup R^{\prime} \) has Hecke-symmetry 
for any set of hyperbolic points \( R \).

In the case of the modular group \( G_{3} \),
hyperbolic fixed points are elements 
of \( \Z\left[ \sqrt{D} \right] \)
for \( D = (a+d)^{2}-4 \)
and
Hecke-conjugation reduces to algebraic conjugation.
For other groups,
a Hecke-conjugate \( \alpha^{\prime} \) is \emph{one} of the algebraic
conjugates of \( \alpha \).


\subsection{Conjugacy class generators}

Schmidt and Sheingorn \cite{MR1362251}
observed that
the non-elliptic conjugacy classes 
in \( G_{p} \)
have representatives that 
are products of
the generators
\( V_{j} = U^{j-1}S \) 
for \( 1 \leq j \leq p-1 \).
In \cite{MR3078226}
we show that every non-elliptic element in \( G_{p} \)
is conjugate to a product 
of conjugacy class generators \( V_{j} \)
that
is unique up to cyclic permutation.
We define the \emph{block length}
of a product of generators
\( W = V_{j_{1}}V_{j_{2}} \cdots V_{j_{t}} \)
to be \( t \),
the length of that product.
Every conjugacy class in \( G_{p} \)
has an associated block length.

For a given \( p \geq 3 \) we have the generators
\begin{align*}
    V_{1} & = S = \left( \begin{matrix} 1 & \lambda \\
	    0 & 1 \end{matrix} \right), \\
    V_{2} & = US = \left( \begin{matrix} \lambda & \lambda^{2}-1 \\
	    1 & \lambda \end{matrix} \right), \\
    V_{3} & = U^{2}S 
	    = \left( \begin{matrix} \lambda^{2}-1 & \lambda^{3}-2\lambda \\
	    \lambda & \lambda^{2}-1 \end{matrix} \right), \\
    & \vdots   \\
    V_{p-3} & = U^{p-4}S  
	    = \left( \begin{matrix} \lambda^{2}-1 & \lambda \\
	    \lambda^{3}-2\lambda & \lambda^{2}-1 \end{matrix} \right), \\
    V_{p-2} & = U^{p-3}S = \left( \begin{matrix} \lambda & 1 \\
	    \lambda^{2}-1 & \lambda \end{matrix} \right), \\
    V_{p-1} & = U^{p-2}S = \left( \begin{matrix} 1 & 0 \\
	    \lambda & 1 \end{matrix} \right).
\end{align*}
Generators
\( V_{1} \) and \( V_{p-1} \)
are parabolic,
as are powers of \( V_{1} \) 
and powers of \( V_{p-1} \).
All other generators and generator products are hyperbolic.


\subsection{\( \lambda \)-continued fractions}

We will use a modification of Rosen's 
continued fractions 
\cite{MR0065632}, 
which are closely associated with the Hecke groups.

For real
\( \alpha \) we put
\( \alpha_{0}=\alpha \) and 
define
\( r_{j} = \left[ \frac{\alpha_{j}}{\lambda} \right] +1 \) and
\( \alpha_{j+1}=\frac{1}{r_{j}\lambda - \alpha_{j}} \)
for \( j \geq 0 \).
Then \( \alpha_{j}=r_{j}\lambda-\frac{1}{\alpha_{j+1}} \) 
for \( j \geq 0 \) 
and 
\( \alpha = r_{0}\lambda-\frac{1}{r_{1}\lambda-\ddots}
          = [r_{0};r_{1},\dots]\)
is the
\emph{\(\lambda_{p} \)-continued fraction} 
(\( \lambda\)-CF)
for \( \alpha \).
An \emph{admissible} \(\lambda\)-CF is one that arises 
from a finite real number by this algorithm.

An admissible 
\( \lambda \)-CF 
has at most \( p-3 \) consecutive ones in any position 
but the beginning, and
has at most \( p-2 \) consecutive ones at the beginning
\cite{MR2582982}.
Schmidt and Sheingorn
\cite{MR1362251}
show that a real number is a non-elliptic fixed point of 
\( G_{p} \) 
if and only if it has a periodic 
\( \lambda \)-CF;
the number is parabolic 
if its \( \lambda \)-CF has period 
\( [\overline{2,\underbrace{1,\ldots,1}_{p-3}}] \),
and hyperbolic
if its \( \lambda \)-CF 
has any other period.
Thus every hyperbolic number \( \alpha \)
has a \( \lambda \)-CF expansion of the form
\begin{equation*}
\alpha 
= [r_{0};r_{1},\ldots,r_{n},\overline{r_{n+1},\ldots,r_{n+m}}],
\end{equation*}
with a period that is not a cyclic permutation
of \( [\overline{2,\underbrace{1,\ldots,1}_{p-3}}] \).
Two hyperbolic numbers are \( G_{p} \)-equivalent
if and only if
they have the same \( \lambda \)-CF period.


\subsection{\( \Z[\lambda] \)-binary quadratic forms}

We let 
\( \mathcal{Q}_{p,D} \) denote the set
of binary quadratic forms
\[ Q(x,y) = Ax^{2}+Bxy+Cy^{2}, \]
with coefficients in  
\( \Z[\lambda_{p}] \)
and discriminant \( D \).
We also denote a form by 
\( Q=[A,B,C] \) and refer to it as a
\( \lambda \)-BQF.
We restrict our attention to hyperbolic forms,
which are
indefinite forms
associated with hyperbolic elements of \( G_{p} \)
as described below.

A Hecke group 
acts on \( \mathcal{Q}_{p,D} \) by
\( \left(Q \circ M \right)(x,y) = Q(a x + b y, c x + d y) \)
for 
\( Q \in \mathcal{Q}_{p,D} \) 
and
\( M = \bigl( \begin{smallmatrix} a &b \\ c & d 
	      \end{smallmatrix} \bigr) \in G_{p} \).
This action 
preserves the discriminant
and
induces an equivalence relation on \( \mathcal{Q}_{p,D} \).
We say that 
\( Q_{1} \) and \( Q_{2} \) are 
\( G_{p} \)-\emph{equivalent}
if there exists 
\( M \in G_{p} \) such that
\( Q_{2} = Q_{1}\circ M \).


\subsection{Matrices, forms, and fixed points}
\label{BigCorrespondence}

We will use a correspondence
between
primitive hyperbolic elements of \( G_{p} \),
certain indefinite \( \lambda \)-BQFs,
and hyperbolic fixed points.
Because the Hecke groups are projective
we may assume that all matrices have positive trace.
Because all positive powers of a matrix
have the same attracting and repelling fixed points
we restrict our attention to primitive matrices.
We use subscripts to indicate connections
between matrices, forms and numbers.
The formal details 
of the isomorphisms appear in \cite{MR2582982}.

In this correspondence 
a primitive hyperbolic matrix 
with positive trace
\[ M = \begin{pmatrix} a & b \\ c & d \end{pmatrix}, \]
corresponds to the 
indefinite \( \lambda \)-BQF
\begin{equation*}
Q_{M} = [A,B,C] = [c,d-a,-b],
\end{equation*}
with discriminant \( D = (a+d)^{2}-4 \).
Because of its connection to \( M \),
we say that such a \( \lambda \)-BQF \( Q_{M} \) is \emph{hyperbolic}.
The form \( Q_{M} \) in turn corresponds to
\begin{equation*}
\alpha_{Q_{M}} = \alpha_{M} = \frac{-B+\sqrt{D}}{2A} = \frac{a-d+\sqrt{D}}{2c},
\end{equation*}
which is the attracting fixed point
of \( M \).
We complete the correspondence
by using the fact that
a hyperbolic fixed point \( \alpha \)
has a \( \lambda \)-CF expansion of the form
\begin{equation*}
\alpha = [r_{0};r_{1},\ldots,r_{n},\overline{r_{n+1},\ldots,r_{n+m}}].
\end{equation*}
If we put
\( V=S^{r_{0}}TS^{r_{1}}T \cdots S^{r_{n}}T \) and
\( W=S^{r_{n+1}}T \cdots S^{r_{n+m}}T \)
the element of \( G_{p} \)
that corresponds to \( \alpha \)
is 
\begin{equation*}
M_{\alpha} = VWV^{-1}.
\end{equation*}
The matrix 
\( M_{\alpha} \) 
is primitive, hyperbolic,
and has
\( \alpha \) as an attracting fixed point.
The \( \lambda \)-BQF
that corresponds to \( \alpha \)
is \( Q_{M_{\alpha}} = Q_{\alpha} \).


\subsection{Hecke-conjugation}

Given a hyperbolic matrix \( M \),
its inverse \( M^{-1} \)
is also hyperbolic and
has the same fixed points,
but with attracting and repelling points reversed.
As a result
the Hecke-conjugate of \( \alpha_{M} \)
is 
\begin{equation*}
\alpha_{M}^{\prime} = \alpha_{M^{-1}}.
\end{equation*}
A calculation shows that
for any hyperbolic \( \alpha \) 
and \( V \in G_{p} \)
we have \( (V\alpha)^{\prime} = V\alpha^{\prime} \).

Given a hyperbolic \( \lambda \)-BQF \( Q = [A,B,C] \)
its negative \( -Q = [-A,-B,-C] \)
is also hyperbolic.
A simple calculation shows that
if \( Q = Q_{M} \) then \( -Q = Q_{M^{-1}} \).
Thus 
the Hecke-conjugate of \( \alpha_{Q} \)
is 
\begin{equation*}
\alpha_{Q}^{\prime} = \alpha_{-Q}.
\end{equation*}


\subsection{Reduction and simplicity}

The positive poles of every rational period function
are ``simple'' numbers associated with simple \( \lambda \)-BQFs,
which are closely related to reduced \( \lambda \)-BQFs.

Zagier's theory of reduction
for classical binary quadratic forms
in \cite{MR0631688}
uses negative classical continued fractions.
We will use a generalization 
to \( \lambda \)-BQFs 
developed in \cite{MR2582982}.

We say that a real number \( \beta \) 
is \emph{\( G_{p} \)-reduced} 
if its \(\lambda_{p} \)-CF expansion 
is purely periodic with a period 
that is not 
\( [\overline{2,\underbrace{1,\ldots,1}_{p-3}}] \).
If \( \beta \) is \( G_{p} \)-reduced
it is hyperbolic, and
we also say that the associated hyperbolic 
\(\lambda \)-BQF \( Q_{\beta} \) is
\( G_{p} \)-\emph{reduced}.

Any hyperbolic \(\lambda \)-BQF \( Q \)
can be transformed into a reduced form
by the action of finitely many elements of \( G_{p} \).
This process maps reduced forms to reduced forms
and so
produces a cycle of reduced forms
in the same equivalence class as \( Q \).
For the associated hyperbolic number \( \beta_{Q} \),
this reduction 
removes the pre-period
then cyclically permutes the period
to produce a cycle of reduced numbers
in the same equivalence class as \( \beta_{Q} \).
Every equivalence class of hyperbolic numbers
has a unique cycle of reduced numbers.

We may characterize reduced numbers
in terms of their size alone.
We show in \cite{MR2582982}
that a hyperbolic fixed point 
\( \beta \) is \( G_{p} \)-reduced
with \( j \) leading ones in its \( \lambda \)-CF
if and only if 
\begin{equation}
	0 < \beta^{\prime} < U^{j+2}(0) < \beta < U^{j+1}(0),
    \label{eq:ReducedNumbers}
\end{equation}
for 
\( 0 \leq j \leq p-3 \).
If the \( \lambda \)-CF 
for a reduced number \( \beta \)
does not have a leading \( 1 \),
then from \eqref{eq:ReducedNumbers} 
we have \( \beta > \lambda = U^{2}(0) \).
If the \( \lambda \)-CF 
for \( \beta \)
has a leading \( 1 \)
we have \( \beta < \lambda \).

We say that a hyperbolic \( \lambda \)-BQF \( Q = [A,B,C] \) 
is \( G_{p} \)-\emph{simple} if
\( A > 0 > C \);
we also say that 
the associated hyperbolic number
\( \alpha_{Q} \) is \( G_{p} \)-\emph{simple}.
A hyperbolic number
\( \alpha \) is \( G_{p} \)-simple 
if and only if 
\( \alpha^{\prime}<0<\alpha \).
Every reduced number \( \beta \)
that is greater than \( \lambda \)
(so its \( \lambda \)-CF does not have a leading \( 1 \))
is associated with at least one simple number
because \( 0 < \beta^{\prime} < \lambda < \beta \)
implies \( (S^{-1}\beta)^{\prime} < 0 < S^{-1}\beta \).

If \( \mathcal{A} \) is an equivalence class 
of hyperbolic \( \lambda \)-BQFs,
we define the corresponding set
of simple numbers
\( \mathcal{Z}_{\mathcal{A}} 
    = \left\{\alpha : Q_{\alpha} \in \mathcal{A}, 
    \alpha\ G_{p}\textrm{-simple} \right\} \).
These sets are nonempty
because every hyperbolic equivalence class 
has a cycle of reduced numbers,
at least one of which 
has a \( \lambda \)-CF with a leading entry 
greater than \( 1 \).
In \cite{MR2582982} we show that
the set of simple numbers 
for a hyperbolic \( \lambda \)-BQF equivalence class  \( \mathcal{A} \)
is 
\begin{equation}
\label{eq:SimpleNumbers}
\mathcal{Z}_{\mathcal{A}}  
= \left\{S^{-i}\beta:
Q_{\beta} \in \mathcal{A} \text{ is \( G_{p} \)-reduced, } 
1 \leq i \leq \left[\frac{\beta}{\lambda} \right]\right\}.
\end{equation}


\subsection{Rational period functions}

For the matrix
\( M = \left( \begin{smallmatrix} * & * \\ c & d \end{smallmatrix} \right) \)
and the function \( f(z) \),
we define the
\emph{weight \( 2k \) slash operator}
\( f \mid_{2k} M = f \mid M \) by
\begin{equation*}
\left( f \mid M \right)(z) = (cz+d)^{-2k}f(Mz).
\end{equation*}
For a fixed \( p \geq 3 \) and positive integer \( k \)
we define a
\emph{rational period function (RPF) of weight \( 2k \) for \( G_{p} \)} 
to be a rational function that satisfies the relations
\begin{equation}
\label{eq:first_relation}
q + q \mid T = 0,
\end{equation}
and
\begin{equation}
\label{eq:second_relation}
q + q \mid U + \cdots + q \mid U^{p-1} = 0.
\end{equation}
This definition
is equivalent to Marvin Knopp's original definition
of rational period functions
for automorphic integrals
\cite{MR0344454}.
The set of rational functions
that satisfy \eqref{eq:first_relation} and \eqref{eq:second_relation}
forms a vector space.

Following Hawkins' insight for RPFs on the modular group \cite{Hawkins_manuscript}
we define an
\emph{irreducible system of poles (ISP)}
to be the minimal set of nonzero poles
forced to occur together
by the relations 
\eqref{eq:first_relation} and \eqref{eq:second_relation}.
For some weights
an RPF with poles in a given ISP
must also have a pole at \( 0 \).

The poles of an RPF on \( G_{p} \) are all real,
and the nonzero poles are all 
hyperbolic fixed points of \( G_{p} \).
The set of positive poles 
in any given ISP
is \( \mathcal{Z}_{\mathcal{A}} \)
for some equivalence class \( \mathcal{A} \)
of hyperbolic forms.

If \( q \) is an RPF of weight \( 2k \) on \( G_{p} \) 
with a pole \emph{only} at zero,
then 
\( q \) must have the form
\cite{MR623004}
\begin{equation}
  q_{k,0}(z) = 
    \begin{cases}
      a_{0}(1-z^{-2k}), & \text{if } 2k \neq 2, \\
      a_{0}(1-z^{-2}) + b_{1}z^{-1}, & \text{if } 2k = 2.
    \end{cases}
    \label{eq:RPFPoleatZero}
\end{equation}
For a nonzero pole
\( \alpha \)
we define
\begin{equation}
q_{k,\alpha}(z)
= PP_{\alpha} 
  \left[ \frac{D^{k/2}}{Q_{\alpha}(z,1)^{k}} \right]
= PP_{\alpha} 
         \left[ \frac{(\alpha-\alpha^{\prime})^{k}}                                       
	 {(z-\alpha)^{k}(z-\alpha^{\prime})^{k}}
                 \right],
    \label{eq:PPatalpha}
\end{equation}
where \(D\) is the discriminant of 
the corresponding \( \lambda \)-BQF \( Q_{\alpha} \).
Using this notation we have
the following expression
for any RPF on \( G_{p} \)
\cite{MR1876701}.
\begin{theorem}
\label{thm:AnyRPFdiscription}
Fix \( p \geq 3 \) and let 
\( \lambda = \lambda_{p} \).
An RPF of weight \( 2k \in 2\Z^{+} \) on \( G_{p} \) is of the form
\begin{equation}
	\label{eq:RPFform}
    q(z) = \sum_{\ell=1}^{L} C_{\ell} 
	    \left(\sum_{\alpha \in Z_{\mathcal{A}_{\ell}}}
	            q_{k,\alpha}(z) 
		- \sum_{\alpha \in Z_{-\mathcal{A}_{\ell}}}
		    q_{k,\alpha^{\prime}}(z)
	    \right) 
	    + c_{0}q_{k,0}(z)
	    + \sum_{n=1}^{2k-1}\frac{c_{n}}{z^{n}},
\end{equation}
where each
\( \mathcal{A}_{\ell} \) is a \( G_{p} \)-equivalence class of 
    \( \lambda \)-BQFs,
\( Z_{\mathcal{A}_{\ell}} \) is the set of positive poles 
    associated with \( \mathcal{A}_{\ell} \),
\( q_{k,\alpha} \) is given by \eqref{eq:PPatalpha},
\( q_{k,0} \) is given by \eqref{eq:RPFPoleatZero},
and the \( C_{\ell} \) and \( c_{n} \) are all constants.
\end{theorem}
The last sum in \eqref{eq:RPFform} is required by
the partial fraction decompositions
in \eqref{eq:PPatalpha}.

If \( \mathcal{A} = -\mathcal{A} \)
then the associated ISP
\( P_{\mathcal{A}} \)
has Hecke-symmetry.
If an RPF of weight \( 2k \)
with poles in \( P_{\mathcal{A}} \)
exists for this case 
it has the form 
\begin{equation}
q(z) = \sum_{\alpha \in Z_{\mathcal{A}}} 
    \left( q_{k,\alpha}(z) - q_{k,\alpha^{\prime}}(z) \right)
		+ \sum_{n=1}^{2k-1}\frac{c_{n}}{z^{n}}.
    \label{eq:RPFsymmetrick}
\end{equation}
If \( k \) is odd then
\begin{equation}
q(z) = \sum_{\alpha \in Z_{\mathcal{A}}} Q_{\alpha}(z,1)^{-k},
    \label{eq:RPFsymmetrickodd}
\end{equation}
is an RPF of weight \( 2k \)
with poles in \( P_{\mathcal{A}} \) \cite{MR1876701}.

If \( \mathcal{A} \neq -\mathcal{A} \)
then the associated ISPs
\( P_{\mathcal{A}} \) and \( P_{-\mathcal{A}} \)
do not have Hecke-symmetry.
If an RPF of weight \( 2k \)
with poles in \( P_{\mathcal{A}} \)
exists for this case 
it has the form 
\begin{equation}
q(z) = \sum_{\alpha \in Z_{\mathcal{A}}}
	            q_{k,\alpha}(z) 
		- \sum_{\alpha \in Z_{-\mathcal{A}}}
		    q_{k,\alpha^{\prime}}(z)
		+ \sum_{n=1}^{2k-1}\frac{c_{n}}{z^{n}}.
    \label{eq:RPFnonsymmetrick}
\end{equation}
The union
\( P_{\mathcal{A}} \cup P_{-\mathcal{A}} \)
does have Hecke-symmetry.
Then
\begin{equation}
q(z) = \sum_{\alpha \in Z_{\mathcal{A}}} Q_{\alpha}(z,1)^{-k}
	- (-1)^{k}\sum_{\alpha \in Z_{-\mathcal{A}}} Q_{\alpha}(z,1)^{-k},
\label{eq:RPFquadraticanyk}
\end{equation}
is an RPF of weight \( 2k \)
for any \( k \)
with poles in 
\( P_{\mathcal{A}} \cup P_{-\mathcal{A}} \) \cite{MR3574636}.

For a fixed \( k \) and ISP \( P_{\mathcal{A}} \),
the results quoted above
characterize nontrivial RPFs
for certain values of \( k \)
and certain kinds of ISPs.
In particular,
\eqref{eq:RPFsymmetrickodd} and \eqref{eq:RPFquadraticanyk}
imply the existence of
a nontrivial RPF for
\begin{enumerate}
\renewcommand{\labelenumi}{(\roman{enumi})}
\item 
\( k \) odd and Hecke-symmetric pole set \( P_{\mathcal{A}} \),
\item 
any \( k \)
and pole set \( P_{\mathcal{A}} \cup P_{-\mathcal{A}} \).
\end{enumerate}
Existence is not guaranteed
and
few RPFs have been constructed
for the other cases
\begin{enumerate}
\renewcommand{\labelenumi}{(\roman{enumi})}
\setcounter{enumi}{2}
\item 
\( k \) even and Hecke-symmetric pole set \( P_{\mathcal{A}} \),
\item 
any \( k \)
and non-Hecke-symmetric pole set \( P_{\mathcal{A}} \).
\end{enumerate}
We will write several new examples 
of RPFs in case (iii) and (iv)
in Section \ref{sec:examples}.


\section{Conjugacy classes and irreducible systems of poles}

We have outlined a correspondence between
primitive hyperbolic matrices,
hyperbolic \( \lambda \)-BQFs,
and hyperbolic fixed points.
The next lemma
shows that the
action of \( G_{p} \) on numbers and \( \lambda \)-BQFs
corresponds to conjugation of matrices.
\begin{lemma}			
\label{lem:Numbers-BQFs-Matrices}
Fix \( p \geq 3 \) and let 
\( \lambda = \lambda_{p} \).
Suppose that \( \alpha \) and \( \beta \)
are hyperbolic numbers
associated with primitive hyperbolic matrices
\( M_{\alpha}, M_{\beta} \in G_{p} \)
and with hyperbolic \( \lambda \)-BQFs
\( Q_{\alpha} \) and \( Q_{\beta} \).
Then for any \( V \in G_{p} \) 
the following statements are equivalent.
\begin{enumerate}
\renewcommand{\labelenumi}{(\alph{enumi})}
\item 
\( \beta = V\alpha \), and
\item 
\( Q_{\beta} = Q_{\alpha}\circ V^{-1} \).
\item 
\( M_{\beta} = VM_{\alpha}V^{-1} \)
\end{enumerate}
\end{lemma}
\begin{proof}
Lemma 7 in \cite{MR2582982}
showed that (a) and (b) are equivalent.
In order to prove that (a) implies (c) 
we use the fact that if (a) holds
then \( \alpha \) and \( \beta \) 
must have the same \( \lambda \)-CF period,
and calculate the matrices 
\( M_{\beta} \)
and 
\( VM_{\alpha}V^{-1} \).
The proof of the converse is another
direct calculation.
\end{proof}

We will use \( \sim \)
to denote equivalence in all three settings:
equivalence of numbers with respect to \( G_{p} \),
equivalence of \( \lambda \)-BQFs with respect to \( G_{p} \),
and conjugacy of elements of \( G_{p} \).
So Lemma \ref{lem:Numbers-BQFs-Matrices} means that 
\( \alpha \sim \beta \) (hyperbolic numbers)
if and only if
\( Q_{\alpha} \sim Q_{\beta} \) (quadratic forms)
if and only if 
\( M_{\alpha} \sim M_{\beta} \) (conjugation).

The correspondence in Section \ref{BigCorrespondence}
holds only for \emph{primitive} hyperbolic elements of \( G_{p} \).
Matrices in a conjugacy class
are either all primitive or all non-primitive,
so
we will describe conjugacy classes themselves
as primitive or not.
Lemma \ref{lem:Numbers-BQFs-Matrices}
has a corollary for equivalence classes.
\begin{corollary}
\label{cor:Numbers-BQFs-Matrices}
Fix \( p \geq 3 \) and let 
\( \lambda = \lambda_{p} \).
The following sets
are in one-to-one correspondence
for the Hecke group \( G_{p} \):
\begin{itemize}
\item 
equivalence classes of hyperbolic numbers,
\item 
equivalence classes of hyperbolic \( \lambda \)-BQFs,
and
\item 
primitive hyperbolic conjugacy classes.
\end{itemize}
\end{corollary}


\subsection{Connecting ISPs and conjugacy classes}
\label{sec:constructISP}

A primitive hyperbolic conjugacy class in \( G_{p} \)
corresponds to a unique 
\( G_{p} \)-equivalence class 
of \( \lambda \)-BQFs \( \mathcal{A} \), 
which in turn is associated 
with 
a unique rational period function ISP \( P_{\mathcal{A}} \).
The elements of
\( P_{\mathcal{A}} \)
all have the same \( \lambda \)-CF period
and thus lie in the same \( G_{p} \)-equivalence class of numbers.
We reverse this 
to find the unique conjugacy classes 
for a given ISP.

We can make this explicit.
Given a primitive hyperbolic conjugacy class in \( G_{p} \)
we write its product of conjugacy class generators.
We pay particular attention to \( V_{1} = S \)
because it plays a special role in the correspondence.
We cyclically permute the generators (by conjugation)  
if necessary 
to write a generator product in the form
\begin{equation}
\label{eq:GenProduct}
W = V_{1}^{m_{1}}V_{j_{1}}V_{1}^{m_{2}}V_{j_{2}} 
	    \cdots V_{1}^{m_{\ell}}V_{j_{\ell}},
\end{equation}
where \( m_{t} \geq 0 \) 
and \( 2 \leq j_{t} \leq p-1 \) for \( 1 \leq t \leq \ell \).
This
\( W \) has \( \ell \) generators that are not \( V_{1} \)
and
\( m_{t} \) copies of \( V_{1} \) preceding each \( V_{j_{t}} \).
Then
\( W \) is conjugate to
\begin{equation*}
M = SWS^{-1}
= S^{m_{1}+2}T(ST)^{j_{1}-2}S^{m_{2}+2}T(ST)^{j_{2}-2} 
	    \cdots S^{m_{\ell}+2}T(ST)^{j_{\ell}-2},
\end{equation*}
which corresponds to the reduced fixed point
\begin{equation}
\label{eq:reducedCF}
\beta_{1}
= [\overline{m_{1}+2;\underbrace{1,1,\ldots, 1}_{j_{1}-2}, 
	m_{2}+2, \underbrace{1,1,\ldots, 1}_{j_{2}-2}, \ldots ,
	m_{\ell}+2, \underbrace{1,1,\ldots, 1}_{j_{\ell}-2} }].
\end{equation}
We observe that 
the \( \lambda \)-CF is admissible
because \( j_{t}-2 \leq p-3 \) for each \( t \)
and that its period is minimal
because the original matrix \( W \) is primitive.
We let \( \mathcal{A} \)
denote the \( G_{p} \)-equivalence class of \( \lambda \)-BQFs
that corresponds to  
the \( G_{p} \)-equivalence class of numbers
\( [\beta_{1}] \).

By \eqref{eq:SimpleNumbers},
the simple numbers in \( [\beta_{1}] \)
are images 
under powers of \( S^{-1} \)
of reduced numbers 
that have \( \lambda \)-CF leading entry greater than \( 1 \).
These are the \( \ell \) reduced numbers 
\begin{equation*}
\beta_{t}
= [\overline{m_{t}+2;\underbrace{1,1,\ldots, 1}_{j_{t}-2}, \ldots
	}],
\end{equation*}
for \( 1 \leq t \leq \ell \).
The simple number(s) for each \( \beta_{t} \) are
the \( m_{t}+1 \) number(s)
\begin{align*}
\alpha_{i}^{(t)}
& = S^{-i}\beta_{t} 
 = [\overline{m_{t}+2-i;\underbrace{1,1,\ldots, 1}_{j_{t}-2}, \ldots }],
\end{align*}
for \( 1 \leq i \leq \left[ \frac{\beta_{t}}{\lambda} \right] = m_{t}+1 \).
Thus the set of simple numbers 
in \( [\beta_{1}] \) is
\begin{equation}
\label{eq:PosPoles}
\mathcal{Z_{A}}
= \left\{ \alpha_{1}^{(1)}, \ldots, \alpha_{m_{1}+1}^{(1)},
    \alpha_{1}^{(2)}, \ldots, \alpha_{m_{2}+1}^{(2)}, 
    \ldots,
    \alpha_{1}^{(\ell)}, \ldots, \alpha_{m_{\ell}+1}^{(\ell)} \right\}.
\end{equation}
These numbers are the positive poles 
for the ISP 
\( P_{\mathcal{A}} = \mathcal{Z_{A}} \cup T\mathcal{Z_{A}} \).

Conversely,
the positive poles of a given ISP
\( P_{\mathcal{A}} \)
are the simple numbers \( \mathcal{Z_{A}} \).
The equivalence class for \( \mathcal{Z_{A}} \)
contains at least one reduced number \( \beta_{1} \)
of the form \eqref{eq:reducedCF}.
As a result, matrices corresponding 
to poles in \( P_{\mathcal{A}} \)
are all in the same conjugacy class as \( M_{\beta_{1}} \sim W \),
where \( W \) is a conjugacy class generator block.
The product \( W \)
(and so the conjugacy class)
is primitive and hyperbolic
because \( M_{\beta} \) is primitive and hyperbolic.

We summarize our result the following Theorem.
\begin{theorem}
\label{thm:CCtoISP}
Fix \( p \geq 3 \) and let \( \lambda = \lambda_{p} \).
There is a one-to-one correspondence
between
primitive hyperbolic conjugacy classes of \( G_{p} \)
and
ISPs for RPFs on \( G_{p} \).
In particular,
a conjugacy class generator product of the form \eqref{eq:GenProduct}
corresponds to 
a set of simple numbers 
\( \mathcal{Z}_{\mathcal{A}} \)
of the form \eqref{eq:PosPoles}.
\end{theorem}

We simplify subsequent calculations 
by observing that
when we translate
between generator products and \( \lambda \)-CFs
every generator (except \( V_{1} \)) 
corresponds to the sequence of \( \lambda \)-CF entries
listed in Table \ref{GenToCF}.
The occurrence of \( V_{1}^{m} \) corresponds to 
the addition of \( m \) 
to the following \( \lambda \)-CF entry,
making that entry \( 2+m \).
\begin{table}[h]
\caption{Translation between generators and \( \lambda \)-CF periods}
\begin{center}
\begin{tabular}{ccc}
\underline{Generator}	&	\underline{\( \lambda \)-CF entry} \\
\( V_{1}^{m} \)		&	add \( m \) \\
\( V_{2} \)		&	\( [\overline{2}] \) \\
\( V_{3} \)		&	\( [\overline{2,1}] \) \\
\( \vdots \)	&	\( \vdots \) \\
\( V_{p-1} \)		&	\( [\overline{2,\underbrace{1,1,\ldots, 1}_{p-3}}] \) \\
\end{tabular}
\end{center}
\label{GenToCF}
\end{table}


\subsection{Counting poles and ISPs}

One consequence of our construction
is that
the number of positive poles in an ISP
is the same as
the block length 
for the corresponding
conjugacy class.
\begin{corollary}			
\label{cor:numberofpoles}
Fix \( p \geq 3 \) and let \( \lambda = \lambda_{p} \).
A primitive hyperbolic conjugacy class in \( G_{p} \)
has block length \( n \)
if and only if
the corresponding RPF ISP 
has \( n \) positive poles.
\end{corollary}
\begin{proof}
We observe that 
\( \sum_{i=1}^{\ell} (m_{i}+1) \)
is the block length of the product \eqref{eq:GenProduct}
and the cardinality of \( \mathcal{Z_{A}} \)
in \eqref{eq:PosPoles}.
\end{proof}

Given a Hecke group \( G_{p} \)
we let \( B_{p}(n) \)
denote the number 
of ISPs for \( G_{p} \)
that have \( n \) positive poles.
By Corollary \ref{cor:numberofpoles},
\( B_{p}(n) \) is also the number 
of primitive hyperbolic conjugacy classes in \( G_{p} \)
with block length \( n \).
Calculating \( B_{p}(n) \)
is (essentially) the problem
of counting the number
of primitive cyclic words
of length \( n \) on \( q \) letters.

Cyclic words are sometimes called necklaces because
we identify each word with its cyclic permutations,
as we have done with conjugacy class generators
by conjugation.
A cyclic words of length \( n \) has period \( n \),
and possibly sub-period \( d \) for some \( d \mid n \).
A word is \emph{primitive}
if it has no nontrivial sub-period.

Let 
\( P_{q}(n) \) denote
the number of primitive words
of length \( n \) on \( q \) letters,
and
\( C_{q}(n) \) denote
the number of primitive \emph{cyclic} words
of length \( n \) on \( q \) letters.
Then \( C_{q}(n) \)
is the number of primitive necklaces with \( n \) beads
chosen from \( q \) colors.
\cite{MR142480}
The following result is well-known.

\begin{lemma}
\label{lem:countPrimitiveNecklaces}
The number of primitive cyclic words
of length \( n \)
with \( q \) letters
is
\begin{equation*}
C_{q}(n) = \frac{1}{n}\sum_{d \mid n}\mu(d)q^{n/d}.
\end{equation*}
\end{lemma}

We give a short proof inspired 
by a proof of the formula 
for \( P_{q}(n) \) in \cite{MR95091}.

\begin{proof}
There are \( q^{n} \) words of length \( n \)
with \( q \) letters.
Every word has a sub-period that divides \( n \),
so
\begin{equation*}
q^{n} = \sum_{d \mid n} P_{q}(d).
\end{equation*}
By M\"{o}bius inversion
we have
\begin{equation*}
P_{q}(n) = \sum_{d \mid n} \mu\left( \frac{n}{d} \right)q^{d}
	= \sum_{d \mid n} \mu\left( d \right)q^{n/d},
\end{equation*}
where \( \mu(n) \) is the M\"{o}bius function.
If we identify cyclic permutations
every word has \( n \) equivalent words,
so
\begin{equation*}
C_{q}(n) = \frac{1}{n}P_{q}(n)
	= \frac{1}{n}\sum_{d \mid n} \mu\left( d \right)q^{n/d}.
\end{equation*}
\end{proof}

For \( n>1 \)
we have \( B_{p}(n) = C_{p-1}(n) \).
But there are two
parabolic conjugacy class generator products
of block length \( 1 \)
(\( V_{1} \) and \( V_{p-1} \)),
so
\begin{equation*}
B_{p}(1) = C_{p-1}(1) - 2 = p-3.
\end{equation*}

\begin{corollary}
\label{lem:countISPs}
Fix \( p \geq 3 \) and let \( \lambda = \lambda_{p} \).
The number of ISPs with \( n \) positive poles for
the Hecke group \( G_{p} \)
is 
\begin{equation*}
B_{p}(n) 
= \begin{cases}
    p-3, & n=1  \\
    C_{p-1}(n), & n>1
\end{cases}.
\end{equation*}
\end{corollary}

In Table \ref{table:ISPcount}
we list the number of ISPs 
with small numbers of positive poles 
for several Hecke groups.

\begin{table}[htp]
\caption{Number of ISPs in \( G_{p} \) with \( n \) positive poles}
\begin{center}
\begin{tabular}{|c|r|r|r|r|r|}
\hline
\multirow{2}{*}{\( n \)} &
\multicolumn{5}{|c|}{\( p \)} \\
\cline{2-6}
 & 3 & 4 & 5 & 6 & 7 \\
\hline
1 & 0 & 1 & 2 & 3 & 4 \\
2 & 1 & 3 & 6 & 10 & 15 \\
3 & 2 & 8 & 20 & 40 & 70 \\
4 & 3 & 18 & 60 & 150 & 315 \\
5 & 6 & 48 & 204 & 624 & 1554 \\
6 & 9 & 116 & 670 & 2580 & 7735 \\
7 & 18 & 312 & 2340 & 11160 & 39990 \\
8 & 30 & 810 & 8160 & 48750 & 209790 \\
\hline
\end{tabular}
\end{center}
\label{table:ISPcount}
\end{table}%


\subsection{ISPs with few poles}

We can write all ISPs that have \( n \) positive poles
by finding all conjugacy class generator products 
with block length \( n \).
We illustrate this by writing
the ISPs with fewest number of poles 
for several groups.
In Section \ref{sec:examples}
we write RPFs for each of these ISPs.

\begin{example}
\label{ex:ISPfewPoles}
Each of the following ISPs has the smallest possible cardinality
for the given group.

\begin{enumerate}
\item 	
The modular group \( G_{3} \)
has no hyperbolic conjugacy class generators,
so it has no ISPs with \( 2 \) poles.
The single hyperbolic product of length \( 2 \)
in \( G_{3} \)
is \( V_{1}V_{2} \).
The corresponding reduced number 
is \( [\overline{3}] = \frac{3+\sqrt{5}}{2} \)
and
the two simple numbers in that class 
are \( \alpha_{1} = [2;\overline{3}] = \frac{1+\sqrt{5}}{2} \)
and \( \alpha_{2} = [1;\overline{3}] = \frac{-1+\sqrt{5}}{2} \).
The ISP is
\( P_{\mathcal{A}} 
= \mathcal{Z}_{\mathcal{A}} \cup T\mathcal{Z}_{\mathcal{A}}
= \left\{ \frac{1\pm\sqrt{5}}{2}, \frac{-1\pm\sqrt{5}}{2}
    \right\} \).

\item 	
The group \( G_{4} \)
has a single hyperbolic generator \( V_{2} \).
The corresponding reduced number 
is \( [\overline{2}] = 1 + \sqrt{2} \)
and
the simple number in that class is \( \alpha = [1;\overline{2}] = 1 \).
The ISP 
is
\( P_{\mathcal{A}} 
= \mathcal{Z}_{\mathcal{A}} \cup T\mathcal{Z}_{\mathcal{A}}
= \left\{ 1, -1 \right\} \).

\item 	
The group \( G_{5} \)
has two hyperbolic generators \( V_{2} \) and \( V_{3} \).
\begin{enumerate}
\item 	
The generator \( V_{2} \)
corresponds to the reduced number 
\( [\overline{2}] = \lambda + \sqrt{\lambda} \).
The simple number in the class 
is \( \alpha = [1;\overline{2}] = \sqrt{\lambda} \),
and the ISP 
is
\( P_{\mathcal{A}_{1}} 
= \mathcal{Z}_{\mathcal{A}_{1}} \cup T\mathcal{Z}_{\mathcal{A}_{1}}
= \left\{ \sqrt{\lambda}, \frac{-1}{\sqrt{\lambda}} \right\} \).

\item 	
The generator \( V_{3} \)
corresponds to the reduced number 
\( [\overline{2,1}] = \lambda + \frac{1}{\sqrt{\lambda}} \).
The simple number in the class 
is \( \beta = [1;\overline{1,2}] = \frac{1}{\sqrt{\lambda}} \),
and the ISP 
is
\( P_{\mathcal{A}_{2}} 
= \mathcal{Z}_{\mathcal{A}_{2}} \cup T\mathcal{Z}_{\mathcal{A}_{2}}
= \left\{ \frac{1}{\sqrt{\lambda}}, -\sqrt{\lambda} \right\} \).
\end{enumerate}

\item 	
The group \( G_{6} \)
has three hyperbolic generators \( V_{2}, V_{3} \), and \( V_{4} \).

\begin{enumerate}
\item 	
The generator \( V_{2} \)
corresponds to the reduced number 
\( [\overline{2}] = \sqrt{3} + \sqrt{2} \).
The simple number in the class is \( \alpha = [1;\overline{2}] = \sqrt{2} \),
and the ISP 
is
\( P_{\mathcal{A}_{1}}
= \mathcal{Z}_{\mathcal{A}_{1}} \cup T\mathcal{Z}_{\mathcal{A}_{1}}
     = \left\{ \sqrt{2}, \frac{-1}{\sqrt{2}} \right\} \).

\item 	
The generator \( V_{3} \)
corresponds to the reduced number 
\( [\overline{2,1}] = \sqrt{3} + 1 \).
The simple number in the class is \( \alpha = [1;\overline{1,2}] = 1 \),
and the ISP 
is
\( P_{\mathcal{A}_{2}} 
= \mathcal{Z}_{\mathcal{A}_{2}} \cup T\mathcal{Z}_{\mathcal{A}_{2}}
= \left\{ 1, -1 \right\} \).

\item 	
The generator \( V_{4} \)
corresponds to the reduced number 
\( [\overline{2,1,1}] = \sqrt{3} + \frac{1}{\sqrt{2}} \).
The simple number in the class
is \( \alpha = [1;\overline{1,1,2}] = \frac{1}{\sqrt{2}} \),
and the ISP 
is
\( P_{\mathcal{A}_{3}} 
= \mathcal{Z}_{\mathcal{A}_{3}} \cup T\mathcal{Z}_{\mathcal{A}_{3}}
= \left\{ \frac{1}{\sqrt{2}}, -\sqrt{2} \right\} \).
\end{enumerate}
\end{enumerate}
\end{example}


\section{Hecke-conjugation}

In this section we study Hecke-conjugation
and Hecke-symmetry,
which play important roles
in the structure of ISPs and RPFs.
We use the results to give several characterizations
of Hecke-symmetry.
Of particular importance
is a connection between
Hecke-conjugate ISPs
and the transpose of conjugacy class generators.


\subsection{Background}

We first state several facts 
from \cite{MR1876701} and \cite{MR2582982}.

Suppose \( \mathcal{A} \) is an equivalence class
of \( \lambda_{p} \)-BQFs.
Then
\( -\mathcal{A} 
     = \left\{-Q \vert Q \in \mathcal{A}\right\} \)
is another equivalence class of forms, 
not necessarily distinct from \( \mathcal{A} \).
The numbers associated with the forms in 
\( -\mathcal{A} \) 
are the Hecke-conjugates of the numbers associated with the forms in 
\( \mathcal{A} \).
The ISP associated with \( \mathcal{A} \)
is 
\begin{align*}
P_{\mathcal{A}}
& = \mathcal{Z}_{\mathcal{A}} \cup T\mathcal{Z}_{\mathcal{A}} \nonumber \\
& = \mathcal{Z}_{\mathcal{A}} \cup \mathcal{Z}_{-\mathcal{A}}^{\prime}.
\end{align*}
An ISP \( P_{\mathcal{A}} \)
has an even number of poles;
the positive half is in \( \mathcal{Z}_{\mathcal{A}} \)
and the negative half is in
\begin{equation*}
\label{eq:rewriteTZ_A}
T\mathcal{Z}_{\mathcal{A}} = \mathcal{Z}_{-\mathcal{A}}^{\prime}.
\end{equation*}
If we take Hecke-conjugates of an ISP \( P_{\mathcal{A}} \)
for \( G_{p} \)
we get
\begin{align*}
P_{\mathcal{A}}^{\prime}
& = \mathcal{Z}_{\mathcal{A}}^{\prime} \cup \mathcal{Z}_{-\mathcal{A}} \\
& = \mathcal{Z}_{-\mathcal{A}} \cup \mathcal{Z}_{-(-\mathcal{A})}^{\prime}  \\
& = P_{\mathcal{-A}},
\end{align*}
which is another ISP
not necessarily distinct from \( P_{\mathcal{A}} \).
If \( P_{\mathcal{A}} = P_{-\mathcal{A}} \)
then \( P_{\mathcal{A}} \) has Hecke-symmetry.
If \( P_{\mathcal{A}} \neq P_{-\mathcal{A}} \)
then \( P_{\mathcal{A}} \) and \( P_{-\mathcal{A}} \)
are Hecke-conjugate ISPs
and
\( P_{\mathcal{A}} \cup P_{-\mathcal{A}} \)
has Hecke-symmetry.


\subsection{Characterizing Hecke-conjugation}

It will be useful to have
several characterizations
of Hecke-conjugate ISPs.
\begin{lemma}
\label{lem:ConjISPcharicterizations}
Fix \( p \geq 3 \) and let 
\( \lambda = \lambda_{p} \).
Suppose that \( \mathcal{A}_{1} \) and \( \mathcal{A}_{2} \) 
are two hyperbolic equivalence classes
of \( \lambda \)-BQFs
with associated ISPs (respectively)
\( P_{\mathcal{A}_{1}} \) and \( P_{\mathcal{A}_{2}} \).
The following statements are equivalent.
\begin{enumerate}
\renewcommand{\labelenumi}{(\alph{enumi})}
\item 
\( P_{\mathcal{A}_{1}}^{\prime} = P_{\mathcal{A}_{2}} \).
\item 
\( -\mathcal{A}_{1} = \mathcal{A}_{2}  \).
\item 
\( T\mathcal{Z}_{\mathcal{A}_{1}}^{\prime} = \mathcal{Z}_{\mathcal{A}_{2}} \).
\item 
\( \alpha_{1}^{\prime} \sim \alpha_{2} \)
for every \( \alpha_{1} \in \mathcal{Z}_{\mathcal{A}_{1}} \)
and \( \alpha_{2} \in \mathcal{Z}_{\mathcal{A}_{2}} \).
\end{enumerate}
\end{lemma}
\begin{proof}
The equivalence of (a), (b) and (c)
is contained in
Section \ref{eq:rewriteTZ_A}.
We show that (b) and (d) are also equivalent.
\begin{itemize}
\item[(b) \( \Rightarrow \) (d)] 
We suppose that
\( -\mathcal{A}_{1} = \mathcal{A}_{2}  \)
and let \( \alpha_{1} \in \mathcal{Z}_{\mathcal{A}_{1}} \)
and \( \alpha_{2} \in \mathcal{Z}_{\mathcal{A}_{2}} \).
Then (c) also holds, so
\[ T\alpha_{1}^{\prime} \in T\mathcal{Z}_{\mathcal{A}_{1}}^{\prime}
    = \mathcal{Z}_{\mathcal{A}_{2}}, \]
so
\( \alpha_{1}^{\prime} \sim \alpha_{2} \).

\item[(d) \( \Rightarrow \) (b)] 
Suppose that 
\( \alpha_{1}^{\prime} \sim \alpha_{2} \)
for every \( \alpha_{1} \in \mathcal{Z}_{\mathcal{A}_{1}} \)
and \( \alpha_{2} \in \mathcal{Z}_{\mathcal{A}_{2}} \).
Fix two such numbers \( \alpha_{1} \) and \( \alpha_{2} \).
Then by Lemma \ref{lem:Numbers-BQFs-Matrices} 
the corresponding \( \lambda \)-BQFs satisfy
\( Q_{\alpha_{1}^{\prime}} \sim Q_{\alpha_{2}} \).
But \( Q_{\alpha_{1}^{\prime}} \in -\mathcal{A}_{1} \)
and \( Q_{\alpha_{2}} \in \mathcal{A}_{2} \),
and because non-disjoint equivalence classes are identical
we have 
\( -\mathcal{A}_{1} = \mathcal{A}_{2} \).
\end{itemize}
\end{proof}

The next lemma shows that
Hecke-conjugation
is related 
to taking the transpose of corresponding matrices.

\begin{lemma}	
\label{lem:mapping_trace}
Fix \( p \geq 3 \) and let \( \lambda = \lambda_{p} \).
Suppose that
\( M \in G_{p} \)
is hyperbolic,
that
\( Q_{M} \) is the associated \( \lambda \)-BQF,
and that \( \alpha_{M} \) is the associated
fixed point.
Then
\begin{equation*}
Q_{M^{\top}} = - Q_{M} \circ T,
\end{equation*}
and
\begin{equation*}
\alpha_{M^{\top}} = T\alpha_{M}^{\prime}.
\end{equation*}
\end{lemma}

\begin{proof}
If we write
\( M = \begin{pmatrix} a & b \\ c & d \end{pmatrix} \)
then
\( M^{\top} = \begin{pmatrix} a & c \\ b & d \end{pmatrix} \).
For the first part we use
\( Q_{M} = [c,d-a,-b] \)
and 
\( Q_{M^{\top}} = [b,d-a,-c] \)
to calculate that
\begin{equation*}
- Q_{M} \circ T 
    = [-c,a-d,b] \circ \begin{pmatrix} 0 & -1 \\ 1 & 0 \end{pmatrix}
    = [b,d-a,-c]
    = Q_{M^{\top}}.
\end{equation*}

For the second result
we use the facts that
\( \alpha_{M} = \frac{a-d+\sqrt{D}}{2c} \)
and
\( \alpha_{M^{\top}} = \frac{a-d+\sqrt{D}}{2b}, \)
where \( D = (a+d)^{2}-4 \),
to calculate that
\begin{align*}
T\alpha_{M}^{\prime}
& = T \left( \frac{a-d-\sqrt{D}}{2c} \right) 
 = \frac{a-d+\sqrt{D}}{2b} 
 = \alpha_{M^{\top}}.
\end{align*}
\end{proof}

We let \( \mathcal{M}_{\mathcal{A}} \)
denote the set of matrices that corresponds
to \( \mathcal{Z}_{\mathcal{A}} \),
that is,
\[ \mathcal{M}_{\mathcal{A}} 
    = \left\{ M_{\alpha} : \alpha 
    \in \mathcal{Z}_{\mathcal{A}} \right\}. \]
By Corollary \ref{cor:Numbers-BQFs-Matrices}
the elements
of a given \( \mathcal{M}_{\mathcal{A}} \) 
are all in the same conjugacy class.  
We also use the transpose symbol
to denote the set of all transpose matrices,
as in
\[ \mathcal{M}_{\mathcal{A}}^{\top} 
    = \left\{ M^{\top}: M \in \mathcal{M}_{\mathcal{A}} \right\}. \]

The next lemma shows
that the transposes of the matrices
for the positive poles in an ISP
are associated with the positive poles
in the conjugate ISP.
\begin{lemma}	
\label{lem:TransposeToConjugate}
Fix \( p \geq 3 \) and let \( \lambda = \lambda_{p} \).
Suppose that \( \mathcal{A}_{1} \) and \( \mathcal{A}_{2} \) 
are two hyperbolic equivalence classes
of \( \lambda \)-BQFs
with associated ISPs (respectively)
\( P_{\mathcal{A}_{1}} \) and \( P_{\mathcal{A}_{2}} \).
Then 
\( P_{\mathcal{A}_{1}} \) and \( P_{\mathcal{A}_{2}} \)
are Hecke-conjugate ISPs
if and only if
\( \mathcal{M}_{\mathcal{A}_{1}}^{\top} = \mathcal{M}_{\mathcal{A}_{2}} \).
\end{lemma}

\begin{proof}
We first suppose that
\( P_{\mathcal{A}_{1}} \) and \( P_{\mathcal{A}_{2}} \)
are Hecke-conjugate ISPs
and let \( M \in \mathcal{M}_{\mathcal{A}_{1}}^{\top} \).
Then
\( M^{\top} \in \mathcal{M}_{\mathcal{A}_{1}} \)
so 
\( \alpha_{M^{\top}} \in \mathcal{Z}_{\mathcal{A}_{1}} \).
Now \( \alpha_{M^{\top}} = T\alpha_{M}^{\prime} \)
by Lemma \ref{lem:mapping_trace},
so
\( T\alpha_{M}^{\prime} \in \mathcal{Z}_{\mathcal{A}_{1}} \).
Thus 
\( \alpha_{M} \in T\mathcal{Z}_{\mathcal{A}_{1}}^{\prime}
	= \mathcal{Z}_{\mathcal{A}_{2}} \)
by Lemma \ref{lem:ConjISPcharicterizations},
so 
\( M \in \mathcal{M}_{\mathcal{A}_{2}} \).
This shows that 
\( \mathcal{M}_{\mathcal{A}_{1}}^{\top} \subseteq \mathcal{M}_{\mathcal{A}_{2}} \).
For containment in the other direction
we can reverse the steps above.
Containment in both directions
means that
\( \mathcal{M}_{\mathcal{A}_{1}}^{\top} = \mathcal{M}_{\mathcal{A}_{2}} \).

Next
we suppose that 
\( \mathcal{M}_{\mathcal{A}_{1}}^{\top} = \mathcal{M}_{\mathcal{A}_{2}} \).
Let \( \alpha_{1} \in \mathcal{Z}_{\mathcal{A}_{1}} \)
and \( \alpha_{2} \in \mathcal{Z}_{\mathcal{A}_{2}} \).
Put
\( M = M_{\alpha_{1}} \), 
so \( \alpha_{1} = \alpha_{M} \).
Then
\( M^{\top} \in 
	\mathcal{M}_{\mathcal{A}_{1}}^{\top} = \mathcal{M}_{\mathcal{A}_{2}} \).
But
\( \alpha_{M^{\top}} = T\alpha_{M}^{\prime} = T\alpha_{1}^{\prime} \)
by Lemma \ref{lem:mapping_trace},
so \( T\alpha_{1}^{\prime} \in \mathcal{Z}_{\mathcal{A}_{2}} \)
and \( \alpha_{1}^{\prime} \sim \alpha_{2} \).
Then Lemma \ref{lem:ConjISPcharicterizations}
gives us that
\( P_{\mathcal{A}_{1}} \) and \( P_{\mathcal{A}_{2}} \)
are Hecke-conjugate ISPs.
\end{proof}

This allows us to characterize 
conjugate ISPs
using the generators of
their corresponding conjugacy classes.
\begin{lemma}	
\label{lem:ConjugateandGeneratorTranspose}
Fix \( p \geq 3 \) and let \( \lambda = \lambda_{p} \).
Suppose that \( W_{1} \) and \( W_{2} \) 
are products of conjugacy class generators in \( G_{p} \)
with associated ISPs (respectively)
\( P_{\mathcal{A}_{1}} \) and \( P_{\mathcal{A}_{2}} \).
Then 
\( P_{\mathcal{A}_{1}} \) and \( P_{\mathcal{A}_{2}} \)
are Hecke-conjugate ISPs
if and only if
\( W_{1}^{\top} \sim W_{2} \).
\end{lemma}

\begin{proof}
By Theorem \ref{thm:CCtoISP}
we have that 
elements of \( \mathcal{M}_{\mathcal{A}_{1}} \)
are in the same conjugacy class as \( W_{1} \)
and 
elements of \( \mathcal{M}_{\mathcal{A}_{2}} \)
are in the same conjugacy class as \( W_{2} \).
Then the result follows from Lemma \ref{lem:TransposeToConjugate}.
\end{proof}

The transpose of a product of matrices
is the product
(in reverse order)
of the transpose matrices,
so in order to use Lemma \ref{lem:ConjugateandGeneratorTranspose}
we need to calculate the transpose of each generator.
But
the transpose of each conjugacy class generator
is another conjugacy class generator.

\begin{lemma}
\label{lem:generator_transpose}
Fix \( p \geq 3 \) and let \( \lambda = \lambda_{p} \).
The conjugacy class generators satisfy
\( V_{j}^{\top} = V_{p-j} \)
for \( 1 \leq j \leq p-1 \).
\end{lemma}
\begin{proof}
In \cite{MR3078226}
we show that
\( V_{j} = \begin{pmatrix} a_{j} & a_{j+1} \\ 
	a_{j-1} & a_{j} \end{pmatrix} \),
where 
\( a_{j} 
= \frac{\sin \left( j\pi/p \right)}
	{\sin \left( \pi/p \right)} \).
Because
\( \sin \left( (p-j)\pi/p \right) = \sin \left( \pi - j\pi/p \right)
= \sin \left( j\pi/p \right) \)
we have
\( a_{p-j} = a_{j} \),
\( a_{p-j+1} = a_{j-1} \), and
\( a_{p-j-1} = a_{j+1} \).
Thus
\begin{align*}
V_{p-j}
& = \begin{pmatrix} a_{p-j} & a_{p-j+1} \\ 
	a_{p-j-1} & a_{p-j} \end{pmatrix} \\
& = \begin{pmatrix} a_{j} & a_{j-1} \\ 
	a_{j+1} & a_{j} \end{pmatrix}  \\
& = V_{j}^{\top}.
\end{align*}
\end{proof}

If a conjugacy class generator product is
\( W_{1} = V_{j_{1}} V_{j_{2}} \cdots V_{j_{n-1}}V_{j_{n}} \),
the transpose is
\begin{align*}
W_{1}^{\top}
& = V_{j_{n}}^{\top}V_{j_{n-1}}^{\top} \cdots V_{j_{2}}^{\top} V_{j_{1}}^{\top} \\
& = V_{p-j_{n}}V_{p-j_{n-1}} \cdots V_{p-j_{2}} V_{p-j_{1}}, 
\end{align*}
another product of conjugacy class generators.
Two generator products 
in the same conjugacy class 
must be cyclic permutations 
of each other.
Thus \( W_{1}^{\top} \sim W_{2} \)
if and only if
\( W_{1}^{\top} \) is a cyclic permutation of \( W_{2} \).
This gives us 
another way to characterize conjugate ISPs
using conjugacy class generators.
\begin{theorem}	%
\label{thm:ConjugateandGeneratorPermutation}
Fix \( p \geq 3 \) and let \( \lambda = \lambda_{p} \).
Suppose that \( W_{1} \) and \( W_{2} \) 
are products of conjugacy class generators in \( G_{p} \)
with associated ISPs (respectively)
\( P_{\mathcal{A}_{1}} \) and \( P_{\mathcal{A}_{2}} \).
Then 
\( P_{\mathcal{A}_{1}} \) and \( P_{\mathcal{A}_{2}} \)
are Hecke-conjugate ISPs
if and only if
\( W_{1}^{\top} \)
is a cyclic permutation of
\( W_{2} \).
\end{theorem}


\subsection{Hecke-symmetry}

An ISP \( P_{\mathcal{A}} \)
is Hecke-symmetric
if it is its own Hecke-conjugate,
so
our characterizations of Hecke-conjugation
give us
characterizations of Hecke-symmetry.

\begin{corollary}[to Lemma \ref{lem:ConjISPcharicterizations}]	
\label{cor:SymmetricISPcharicterizations}
Fix \( p \geq 3 \) and let 
\( \lambda = \lambda_{p} \).
Suppose that \( \mathcal{A} \) 
is a hyperbolic equivalence class
of \( \lambda \)-BQFs
with associated ISP
\( P_{\mathcal{A}} \).
The following statements are equivalent.
\begin{enumerate}
\renewcommand{\labelenumi}{(\alph{enumi})}
\item 
\( P_{\mathcal{A}}^{\prime} = P_{\mathcal{A}} \).
\item 
\( -\mathcal{A} = \mathcal{A}  \).
\item 
\( T\mathcal{Z}_{\mathcal{A}}^{\prime} = \mathcal{Z}_{\mathcal{A}} \).
\item 
\( \alpha^{\prime} \sim \alpha \)
for every \( \alpha \in \mathcal{Z}_{\mathcal{A}} \).
\end{enumerate}
\end{corollary}

There is a connection between Hecke-symmetry of ISPs
and transpose properties of conjugacy class generators.
The next theorem is a corollary
to 
Lemma \ref{lem:TransposeToConjugate}, 
Lemma \ref{lem:ConjugateandGeneratorTranspose}, and
Theorem \ref{thm:ConjugateandGeneratorPermutation}.

\begin{theorem}	
\label{thm:TransposeAndSymmetry}
Fix \( p \geq 3 \) and let \( \lambda = \lambda_{p} \).
Suppose that \( \mathcal{A} \) 
is a hyperbolic equivalence class
of \( \lambda \)-BQFs
with associated ISP
\( P_{\mathcal{A}} \).
The following statements
are equivalent.

\begin{enumerate}
\renewcommand{\labelenumi}{(\alph{enumi})}
\item 
\( P_{\mathcal{A}} \) has Hecke-symmetry.
\item 
\( \mathcal{M}_{\mathcal{A}}^{\top} = \mathcal{M}_{\mathcal{A}} \).
\item 
\( W^{\top} \sim W \).
\item 
\( W^{\top} \)
is a cyclic permutation of
\( W \).
\end{enumerate}
\end{theorem}

We can easily determine the symmetry properties
of an ISP \( P_{\mathcal{A}} \)
from the associated conjugacy class generator product \( W \).
If \( W^{\top} \)
is a cyclic permutation of \( W \)
then \( P_{\mathcal{A}} \) has Hecke-symmetry.
If \( W^{\top} \)
is not a cyclic permutation of \( W \)
then \( P_{\mathcal{A}} \) does not have Hecke-symmetry,
\( P_{-\mathcal{A}} \) is the ISP 
that is Hecke-conjugate to \( P_{\mathcal{A}} \),
and the union 
\( P_{\mathcal{A}} \cup P_{-\mathcal{A}} 
	= P_{\mathcal{A}} \cup P_{\mathcal{A}}^{\prime} \)
does have Hecke-symmetry.

\begin{example}
\label{ex:ISPsymmetry}
We use conjugacy class generators 
to determine Hecke-symmetry
for the ISPs in
Example \ref{ex:ISPfewPoles}.
\begin{enumerate}
\item 	
In \( G_{3} \)
we found that \( W = V_{1}V_{2} \)
gave us the ISP 
\( P_{\mathcal{A}} 
= \left\{ \frac{1\pm\sqrt{5}}{2}, \frac{-1\pm\sqrt{5}}{2} \right\} \).
We calculate that
\( W^{\top} = V_{2}^{\top}V_{1}^{\top} = V_{1}V_{2} = W \),
so 
\( P_{\mathcal{A}} \)
has Hecke-symmetry.

\item 	
In \( G_{4} \)
we found that \( W = V_{2} \)
gave us the ISP 
\( P_{\mathcal{A}} = \left\{ 1, -1 \right\} \).
We calculate that
\( W^{\top} = V_{2}^{\top} = V_{2} = W \),
so 
\( P_{\mathcal{A}} \)
has Hecke-symmetry.

\item 	
In \( G_{5} \)
the hyperbolic generators 
are \( V_{2} \) and \( V_{3} \).
\begin{enumerate}
\renewcommand{\labelenumi}{(\alph{enumi})}
\item 	
The generator \( W = V_{2} \)
gave us the ISP
\( P_{\mathcal{A}_{1}} 
    = \left\{ \sqrt{\lambda}, \frac{-1}{\sqrt{\lambda}} \right\} \).
We calculate that
\( W^{\top} = V_{2}^{\top} = V_{3} \not\sim W \),
so 
\( P_{\mathcal{A}_{1}} \)
does not have Hecke-symmetry.
  
\item 	
The generator \( W = V_{3} \)
gave us the ISP
\( P_{\mathcal{A}_{2}} 
    = \left\{ \frac{1}{\sqrt{\lambda}}, -\sqrt{\lambda} \right\} \).
We calculate that
\( W^{\top} = V_{3}^{\top} = V_{2} \not\sim W \),
so 
\( P_{\mathcal{A}_{2}} \)
does not have Hecke-symmetry.
\end{enumerate}

The fact that \( V_{2}^{\top} = V_{3} \)
means that ISPs 
\( P_{\mathcal{A}_{1}} \) and \( P_{\mathcal{A}_{2}} \)
are Hecke-conjugate
and their union has Hecke-symmetry.

\item 	
In \( G_{6} \)
the hyperbolic generators 
are \( V_{2}, V_{3} \) and \( V_{4} \).
\begin{enumerate}
\renewcommand{\labelenumi}{(\alph{enumi})}
\item 	
The generator \( W = V_{2} \)
gave us the ISP
\( P_{\mathcal{A}_{1}}
     = \left\{ \sqrt{2}, \frac{-1}{\sqrt{2}} \right\} \).
We calculate that
\( W^{\top} = V_{2}^{\top} = V_{4} \not\sim W \),
so the corresponding ISP 
\( P_{\mathcal{A}_{1}}
     = \left\{ \sqrt{2}, \frac{-1}{\sqrt{2}} \right\} \)
does not have Hecke-symmetry.

\item 	
The generator \( W = V_{3} \)
gave us the ISP
\( P_{\mathcal{A}_{2}} = \left\{ 1, -1 \right\} \).
We calculate that
\( W^{\top} = V_{3}^{\top} = V_{3} = W \),
so 
\( P_{\mathcal{A}_{2}} \)
has Hecke-symmetry.

\item 	
The generator \( W = V_{4} \)
gave us the ISP
\( P_{\mathcal{A}_{3}} = \left\{ \frac{1}{\sqrt{2}}, -\sqrt{2} \right\} \).
We calculate that
\( W^{\top} = V_{4}^{\top} = V_{2} \not\sim W \),
so 
\( P_{\mathcal{A}_{3}} \)
does not have Hecke-symmetry.
\end{enumerate}

The fact that \( V_{2}^{\top} = V_{4} \)
means that ISPs 
\( P_{\mathcal{A}_{1}} \) and \( P_{\mathcal{A}_{3}} \)
are Hecke-conjugate
and their union has Hecke-symmetry.
\end{enumerate}
\end{example}


\section{Examples of rational period functions}
\label{sec:examples}

Our procedure for finding ISPs
takes us a long way toward constructing 
RPFs for Hecke groups.
In this section
we write several rational period functions
for the ISPs we've already found,
as well as for two new ISPs.
Most of the RPFs are new.

\subsection{Simple constructions}

It is is straightforward to write
rational period functions
for two classes of poles
and weights.
\begin{enumerate}
\renewcommand{\labelenumi}{(\roman{enumi})}
\item 
If \( P_{\mathcal{A}} \) 
has Hecke-symmetry,
an RPF with poles \( P_{\mathcal{A}} \) 
is given by \eqref{eq:RPFsymmetrickodd}
for any odd \( k \).
\item 
If 
\( P_{\mathcal{A}_{1}} \) and \( P_{\mathcal{A}_{2}} \)
non-Hecke-symmetric but conjugate to each other
then 
an RPF with poles \( P_{\mathcal{A}_{1}} \cup P_{\mathcal{A}_{2}} \) 
is given by \eqref{eq:RPFquadraticanyk}
for any \( k \).
\end{enumerate}


\begin{example}	
\label{ex:RPFsFewPoles}
Consider the ISPs
from Examples \ref{ex:ISPfewPoles} and \ref{ex:ISPsymmetry}.
The Hecke-symmetric ISPs
give us RPFs with the smallest
number of poles in \( G_{3} \), \( G_{4} \), and \( G_{6} \),
but only for \( k \) odd.
The non-Hecke-symmetric ISPs
give us RPFs with more poles
for any \( k \)
when combined with Hecke-conjugate ISPs
for \( G_{5} \) and \( G_{6} \).
The RPF for \( G_{3} \) is well-known;
the others are new.
\begin{enumerate}
\item 	
In \( G_{3} \)
we found the Hecke-symmetric ISP 
\( P_{\mathcal{A}} 
= \left\{ \frac{1\pm\sqrt{5}}{2}, \frac{-1\pm\sqrt{5}}{2} \right\} \).
If \( k \) is odd
an RPF of weight \( 2k \) 
and poles \( P_{\mathcal{A}} \)
is 
\begin{equation}
\label{eq:MKsRPF}
q(z) = \frac{1}{(z^{2}-z-1)^{k}} + \frac{1}{(z^{2}+z-1)^{k}}.
\end{equation}
This was the first known RPF 
with nonzero poles
for the modular group,
constructed by Marvin Knopp
\cite{MR0485700}.

\item 	
\label{ex3part2}
In \( G_{4} \)
we found the Hecke-symmetric ISP 
\( P_{\mathcal{A}} = \left\{ 1, -1 \right\} \).
If \( k \) is odd
an RPF of weight \( 2k \) 
and poles \( P_{\mathcal{A}} \)
is 
\begin{equation}
\label{eq:RPF-G4easy}
q(z) = \frac{1}{(z^{2}-1)^{k}}.
\end{equation}

\item 	
In \( G_{5} \)
we found the ISPs
\( P_{\mathcal{A}_{1}} 
    = \left\{ \sqrt{\lambda}, \frac{-1}{\sqrt{\lambda}} \right\} \)
and
\( P_{\mathcal{A}_{2}} 
    = \left\{ \frac{1}{\sqrt{\lambda}}, -\sqrt{\lambda} \right\} \).
Neither ISP has Hecke-symmetry
but they are conjugate to each other.
A rational period function 
of weight \( 2k \) 
(for any \( k \))
with poles in
\( P_{\mathcal{A}_{1}} \cup P_{\mathcal{A}_{2}} \)
is
\begin{equation}
\label{ex:RPF-G5easy}
q(z) = \frac{1}{(z^{2}-\lambda)^{k}} -\frac{(-1)^{k}}{(\lambda z^{2}-1)^{k}}.
\end{equation}

\item 	
\label{ex3part4}
In \( G_{6} \)
we found the ISPs
\( P_{\mathcal{A}_{1}}
     = \left\{ \sqrt{2}, \frac{-1}{\sqrt{2}} \right\} \),
\( P_{\mathcal{A}_{2}} = \left\{ 1, -1 \right\} \),
and
\( P_{\mathcal{A}_{3}} = \left\{ \frac{1}{\sqrt{2}}, -\sqrt{2} \right\} \).

The ISP \( P_{\mathcal{A}_{2}} \)
has Hecke-symmetry,
so if \( k \) is odd
an RPF of weight \( 2k \) 
and poles \( P_{\mathcal{A}_{2}} \)
is 
\begin{equation}
\label{eq:RPF-G6easy}
q(z) = \frac{1}{(z^{2}-1)^{k}}.
\end{equation}

The ISPs
\( P_{\mathcal{A}_{1}} \) and \( P_{\mathcal{A}_{3}} \)
do not have Hecke-symmetry
but they are conjugate to each other.
A rational period function 
of weight \( 2k \) 
(for any \( k \))
with poles in
\( P_{\mathcal{A}_{1}} \cup P_{\mathcal{A}_{3}} \)
is
\begin{equation}
\label{ex:RPF-G6easy-nonHS}
q(z) = \frac{1}{(z^{2}-2)^{k}} -\frac{(-1)^{k}}{(2 z^{2}-1)^{k}}.
\end{equation}
\end{enumerate}
\end{example}


We could start with \emph{any} conjugacy class generator block,
write the corresponding ISP,
then use 
\eqref{eq:RPFsymmetrickodd} or \eqref{eq:RPFquadraticanyk}
to write an RPF.
We illustrate this in the next example.

\begin{example}		
We start with two similar generator products
of block length \( 3 \)
in \( G_{6} \)
and write their corresponding ISPs
corresponding RPFs.

\begin{enumerate}
\renewcommand{\labelenumi}{(\alph{enumi})}
\item 	
The generator block \( W = V_{1}V_{3}V_{5} \)
satisfies
\( W^{\top} = (V_{1}V_{3}V_{5})^{\top} = V_{1}V_{3}V_{5} = W \),
so the corresponding ISP 
has Hecke-symmetry
and the corresponding RPF has \( 3 \) positive poles.
The reduced number for
\( W = V_{1}V_{3}V_{5} \)
is \( [\overline{3,1,2,1,1,1}] \),	
so
the simple numbers in the class are 
\( [2;\overline{1,2,1,1,1,3}] \),
\( [1;\overline{1,2,1,1,1,3}] \),
and
\( [1;\overline{1,1,1,3,1,2}] \).
The ISP 
is
\begin{align*}
P_{\mathcal{A}} 
& = \mathcal{Z}_{\mathcal{A}} \cup T\mathcal{Z}_{\mathcal{A}} \\
& = \left\{ \frac{2+\sqrt{7}}{\sqrt{3}},
	\frac{-1+\sqrt{7}}{\sqrt{3}}, \frac{-1+\sqrt{7}}{2\sqrt{3}} \right\}
	\bigcup \left\{ \frac{2-\sqrt{7}}{\sqrt{3}},
	\frac{-1-\sqrt{7}}{2\sqrt{3}}, \frac{-1-\sqrt{7}}{\sqrt{3}} \right\}.
\end{align*}
If \( k \) is odd
an RPF on \( G_{6} \)
of weight \( 2k \)
and poles \( P_{\mathcal{A}_{1}} \cup P_{\mathcal{A}_{2}} \)
is\begin{equation*}
q_{1}(z) = \frac{1}{(3\sqrt{3}z^{2}-12z-3\sqrt{3})^{k}}
    + \frac{1}{(3\sqrt{3}z^{2}+6z-6\sqrt{3})^{k}} 
    + \frac{1}{(6\sqrt{3}z^{2}+6z-3\sqrt{3})^{k}},
\end{equation*}
which we simplify to the RPF
\begin{align*}
q(z)
& = 3^{k}q_{1}(z) \\
& = \frac{1}{(\sqrt{3}z^{2}-4z-\sqrt{3})^{k}}
    + \frac{1}{(\sqrt{3}z^{2}+2z-2\sqrt{3})^{k}} 
    + \frac{1}{(2\sqrt{3}z^{2}+2z-\sqrt{3})^{k}}.
\end{align*}

\item 	
The generator block \( W_{1} = V_{1}V_{2}V_{5} \)
satisfies
\( W^{\top} = (V_{1}V_{2}V_{5})^{\top} = V_{1}V_{4}V_{5} = W_{2} \).
Now \( W_{1} \not\sim W_{2} \)
so the conjugate ISPs 
\( P_{\mathcal{A}_{1}} \) and \( P_{\mathcal{A}_{2}} \)
do not have Hecke-symmetry.
Each ISP has \( 3 \) positive poles
so the union has \( 6 \) positive poles.

The reduced number for
\( W_{1} = V_{1}V_{2}V_{5} \)
is \( [\overline{3,2,1,1,1}] \), 
so the 
simple numbers in \( P_{\mathcal{A}_{1}}  \)
are
\( [2;\overline{2,1,1,1,3}] \),
\( [1;\overline{2,1,1,1,3}] \),
and
\( [1;\overline{1,1,1,3,2}] \).
The ISP for \( W_{1} \) is
\begin{align*}
P_{\mathcal{A}_{1}} 
& = \mathcal{Z}_{\mathcal{A}} \cup T\mathcal{Z}_{\mathcal{A}} \\
& = \left\{ \frac{3 \sqrt{3}+\sqrt{47}}{4},
	\frac{-\sqrt{3}+\sqrt{47}}{4},
	\frac{-2\sqrt{3}+\sqrt{47}}{7} \right\} \\
&  \hspace{20pt} \bigcup
	\left\{ \frac{3\sqrt{3}-\sqrt{47}}{5},
	\frac{-\sqrt{3}-\sqrt{47}}{11},
	\frac{-2 \sqrt{3}-\sqrt{47}}{5} \right\}.
\end{align*}

The reduced number for
\( W_{2} = V_{1}V_{4}V_{5} \)
is \( [\overline{3,1,1,2,1,1,1}] \),
so the 
simple numbers in \( P_{\mathcal{A}_{2}}  \)
are
\( [2;\overline{1,1,2,1,1,1,3}] \),
\( [1;\overline{1,1,2,1,1,1,3}] \),
and
\( [1;\overline{1,1,1,3,1,1,2}] \).
The ISP for \( W_{2} \) is
\begin{align*}
P_{\mathcal{A}_{2}} 
& = \left\{ \frac{3\sqrt{3}+\sqrt{47}}{5},
	\frac{-2 \sqrt{3}+\sqrt{47}}{5},
	\frac{-\sqrt{3}+\sqrt{47}}{11} \right\} \\
&  \hspace{20pt} \bigcup
	\left\{ \frac{3 \sqrt{3}-\sqrt{47}}{4}, 
	\frac{-2\sqrt{3}-\sqrt{47}}{7},
	\frac{-\sqrt{3}-\sqrt{47}}{4} \right\}.
\end{align*}
An RPF on \( G_{6} \)
of weight \( 2k \) (for any \( k \))
with poles \( P_{\mathcal{A}_{1}} \cup P_{\mathcal{A}_{2}} \)
is
\begin{align*}
q(z) 
& = \frac{1}{(4z^{2}-6\sqrt{3}z-5)^{k}}
	+ \frac{1}{(4z^{2}+2\sqrt{3}z-11)^{k}}
    + \frac{1}{(7z^{2}+4\sqrt{3}z-5)^{k}} \\
& \hspace{20pt} - \frac{(-1)^{k}}{(5z^{2}-6\sqrt{3}z-4)^{k}}
	- \frac{(-1)^{k}}{(5z^{2}+4\sqrt{3}z-7)^{k}}
    - \frac{(-1)^{k}}{(11z^{2}+2\sqrt{3}z-4)^{k}}.
\end{align*}
\end{enumerate}
\end{example}


\subsection{More complicated constructions}

It
is is more challenging to write
rational period functions 
with poles in a single ISP \( P_{\mathcal{A}} \)
that satisfies
\begin{enumerate}
\renewcommand{\labelenumi}{(\roman{enumi})}
\setcounter{enumi}{2}
\item 
\( P_{\mathcal{A}} \)
has Hecke-symmetry and \( k \) is even,
or
\item 
\( P_{\mathcal{A}} \)
does not have Hecke-symmetry.
\end{enumerate}
The only previously published examples known to this author are
\begin{itemize}
\item 
several RPFs of weight \( 2 \) or \( 4 \) on \( G_{3} \)
with non-Hecke-symmetric ISPs  \cite{Hawkins_manuscript}, 
\item 
an RPF of weight \( 4 \) on \( G_{3} \)
with a Hecke-symmetric ISP  \cite{Hawkins_manuscript,MR1023922}, and
\item 
an RPF of weight \( 2 \) on \( G_{4} \)
with a non-Hecke-symmetric ISP  \cite{MR3574636}.
\end{itemize}

We offer
new examples for these cases
for the ISPs
in Examples \ref{ex:ISPfewPoles} and \ref{ex:ISPsymmetry}.
We use the fact that
if 
there is an RPF of weight \( 2k \)
with \( k \) even 
and poles in a Hecke-symmetric ISP \( P_{\mathcal{A}} \) 
it must have the given by \eqref{eq:RPFsymmetrick}.
If 
there is an RPF of weight \( 2k \)
and poles in a non-Hecke-symmetric ISP \( P_{\mathcal{A}} \) 
it must have the given by \eqref{eq:RPFnonsymmetrick}.

\begin{example}	
\label{ex:RPFsFewPolesHard}
We again consider the ISPs
from Examples \ref{ex:ISPfewPoles} and \ref{ex:ISPsymmetry}.

\begin{enumerate}
\item 	
In \( G_{3} \)
we found one Hecke-symmetric ISP 
with two positive poles,
\( P_{\mathcal{A}} 
= \left\{ \alpha_{1}, \alpha_{2} \right\} 
\cup \left\{ \alpha_{1}^{\prime}, \alpha_{2}^{\prime} \right\} 
= \left\{ \frac{1\pm\sqrt{5}}{2}, \frac{-1\pm\sqrt{5}}{2} \right\} \).
In Example \ref{ex:RPFsFewPoles}
we wrote an RPF of weight \( 2k \) for any odd \( k \) and this ISP.

Now we write an RPF of weight \( 4 \) (\( k=2 \))
whose ISP is \( P_{\mathcal{A}} \).
By \eqref{eq:RPFsymmetrick}
such an RPF 
must have the form
\begin{align*}
q(z)
& = q_{2,\alpha_{1}}(z) - q_{2,\alpha_{1}^{\prime}}(z)
    + q_{2,\alpha_{2}}(z) - q_{2,\alpha_{2}^{\prime}}(z)
    + \sum_{n=1}^{3}\frac{c_{n}}{z^{n}}.
\end{align*}
We use use \eqref{eq:PPatalpha}
and partial fraction decomposition
to write
\begin{align*}
q(z)
& = \frac{-2/\sqrt{5}}{z-\alpha_{1}} + \frac{1}{(z-\alpha_{1})^{2}} 
	- \frac{2/\sqrt{5}}{z-\alpha_{1}^{\prime}} 
		- \frac{1}{(z-\alpha_{1}^{\prime})^{2}} \\
&	\hspace{30pt}	+ \frac{-2/\sqrt{5}}{z-\alpha_{2}} + \frac{1}{(z-\alpha_{2})^{2}} 
	- \frac{2/\sqrt{5}}{z-\alpha_{2}^{\prime}} 
		- \frac{1}{(z-\alpha_{2}^{\prime})^{2}} 
	+ \sum_{n=1}^{3}\frac{c_{n}}{z^{n}}.
\end{align*}
Substitution into
the two relations
\eqref{eq:first_relation} and \eqref{eq:second_relation}
shows both are satisfied
if \( c_{1} = 8/\sqrt{5} \)
and \( c_{2} = c_{3} = 0 \).
Thus 
\begin{align*}
q(z)
& = \frac{-2/\sqrt{5}}{z-\alpha_{1}} + \frac{1}{(z-\alpha_{1})^{2}} 
	- \frac{2/\sqrt{5}}{z-\alpha_{1}^{\prime}} 
		- \frac{1}{(z-\alpha_{1}^{\prime})^{2}} \\
&	\hspace{30pt}	- \frac{2/\sqrt{5}}{z-\alpha_{2}} + \frac{1}{(z-\alpha_{2})^{2}} 
	- \frac{2/\sqrt{5}}{z-\alpha_{2}^{\prime}} 
		- \frac{1}{(z-\alpha_{2}^{\prime})^{2}} 
	+ \frac{8/\sqrt{5}}{z}, 
\end{align*}
is a rational period function
of weight \( 4 \)
on the modular group \( G_{3} \)
with the Hecke-symmetric ISP
\( P_{\mathcal{A}} 
= \left\{ \frac{1\pm\sqrt{5}}{2}, \frac{-1\pm\sqrt{5}}{2} \right\} \).
Hawkins first wrote this RPF in \cite{Hawkins_manuscript}
and Knopp published it in \cite{MR1023922}.

\item 	
In \( G_{4} \)
we found one Hecke-symmetric ISP 
with one positive pole,
\( P_{\mathcal{A}} = \left\{ 1, -1 \right\} \).
We wrote
an RPF 
of weight \( 2k \) for any odd \( k \)
and this ISP
in Example \ref{ex:RPFsFewPoles}.
An RPF of weight \( 2k \)
for even \( k \) 
must have the form
\begin{align*}
q(z)
& = q_{k,1}(z) - q_{k,-1}(z)
    + \sum_{n=1}^{2k-1}\frac{c_{n}}{z^{n}}.
\end{align*}
We let \( k=2 \),
use \eqref{eq:PPatalpha}
and then partial fractions
to find
that \( q \)
has the form
\begin{equation*}
q(z) = \frac{-1}{z-1} + \frac{1}{(z-1)^{2}}
	- \frac{1}{z-1} - \frac{1}{(z+1)^{2}}
	+ \sum_{n=1}^{3}\frac{c_{n}}{z^{n}}.
\end{equation*}
We substitute and find that the two relations
\eqref{eq:first_relation} and \eqref{eq:second_relation}
are satisfied if
\( c_{1}=2 \) and \( c_{2}=c_{3}=0 \).
Thus
\begin{equation*}
q(z) = \frac{-1}{z-1} + \frac{1}{(z-1)^{2}}
	- \frac{1}{z-1} - \frac{1}{(z+1)^{2}}
	+ \frac{2}{z},
\end{equation*}
is an RPF of weight \( 4 \) on \( G_{4} \)
with Hecke-symmetric ISP 
\( P_{\mathcal{A}} = \left\{ 1, -1 \right\} \).

\item 	
In \( G_{5} \)
we found non-Hecke-symmetric ISPs
with one positive pole each,
\( P_{\mathcal{A}_{1}} 
    = \left\{ \sqrt{\lambda}, \frac{-1}{\sqrt{\lambda}} \right\} \)
and
\( P_{\mathcal{A}_{2}} 
    = \left\{ \frac{1}{\sqrt{\lambda}}, -\sqrt{\lambda} \right\} \).
In Example \ref{ex:RPFsFewPoles}
we wrote 
an RPF of weight \( 2k \) (for any \( k \))
with poles 
in \( P_{\mathcal{A}_{1}}  \cup P_{\mathcal{A}_{2}}  \).

\begin{enumerate}
\renewcommand{\labelenumi}{(\alph{enumi})}
\item 	
By \eqref{eq:RPFnonsymmetrick}
an RPF of weight \( 2k \)
with poles only in \( P_{\mathcal{A}_{1}}  \) 
has the form  
\begin{align*}
q(z)
& = q_{k,\sqrt{\lambda}}(z) - q_{k,-1/\sqrt{\lambda}}(z)
	+ \sum_{n=1}^{2k-1}\frac{c_{n}}{z^{n}}.
\end{align*}
We let \( k=1 \)
and substitute to find that the two relations
\eqref{eq:first_relation} and \eqref{eq:second_relation}
are satisfied if \( c_{1} = 0 \).
This gives us that
\begin{equation}
\label{eq:RPF-G5hard-1}
q(z) = \frac{1}{z-\sqrt{\lambda}} - \frac{\sqrt{\lambda}}{\sqrt{\lambda}z+1}
\end{equation}
is a rational period function of weight \( 2 \)
with the non-Hecke-symmetric ISP \( P_{\mathcal{A}_{1}} \).    
     
\item 	
By \eqref{eq:RPFnonsymmetrick}
an RPF of weight \( 2k \)
with poles only in \( P_{\mathcal{A}_{2}}  \) 
has the form  
\begin{align*}
q(z)
& = q_{k,1/\sqrt{\lambda}}(z) - q_{k,-\sqrt{\lambda}}(z)
		+ \sum_{n=1}^{2k-1}\frac{c_{n}}{z^{n}}.
\end{align*}
We let \( k=1 \)
and substitute to find that the two relations
\eqref{eq:first_relation} and \eqref{eq:second_relation}
are satisfied if \( c_{1} = 0 \).
Thus
\begin{equation}
\label{eq:RPF-G5hard-2}
q(z) = \frac{\sqrt{\lambda}}{\sqrt{\lambda}z-1} - \frac{1}{z+\sqrt{\lambda}}
\end{equation}
is a rational period function of weight \( 2 \)
with the non-Hecke-symmetric ISP \( P_{\mathcal{A}_{2}} \).
\end{enumerate}

\item 	
In \( G_{6} \)
we found three ISPs with one positive pole each.
The ISP
\( P_{\mathcal{A}_{2}}
     = \left\{ 1, -1 \right\} \)
is Hecke-symmetric,
and ISP
\( P_{\mathcal{A}_{1}}
     = \left\{ \sqrt{2}, \frac{-1}{\sqrt{2}} \right\} \)
and
\( P_{\mathcal{A}_{3}} = \left\{ \frac{1}{\sqrt{2}}, -\sqrt{2} \right\} \)
do not have Hecke-symmetry
but are conjugate to each other.
In Example \ref{ex:RPFsFewPoles}
we wrote 
an RPF of weight \( 2k \) for \( k \) odd
with poles 
in \( P_{\mathcal{A}_{2}} \)
and 
an RPF 
of weight \( 2k \) (for any \( k \))
with poles \( P_{\mathcal{A}_{1}} \cup P_{\mathcal{A}_{3}} \).

\begin{enumerate}
\renewcommand{\labelenumi}{(\alph{enumi})}
\item 	
An RPF of weight \( 2k \)
with poles only in \( P_{\mathcal{A}_{1}}  \) 
has the form  
\begin{align*}
q(z)
& = q_{k,\sqrt{2}}(z) - q_{k,-1/\sqrt{2}}(z)
	+ \sum_{n=1}^{2k-1}\frac{c_{n}}{z^{n}}.
\end{align*}
For \( k=1 \)
we have that
\begin{equation}
\label{eq:RPF-G6hard-1}
q(z) = \frac{1}{z-\sqrt{2}} - \frac{\sqrt{2}}{\sqrt{2}z+1},
\end{equation}
is a rational period function of weight \( 2 \)
with the non-Hecke-symmetric ISP \( P_{\mathcal{A}_{1}} \).

\item 	
By \eqref{eq:RPFsymmetrick}
an RPF of weight \( 2k \)
for even \( k \) 
with poles in the Hecke-symmetric ISP \( P_{\mathcal{A}_{2}} \)
must have the form
\begin{align*}
q(z)
& = q_{k,1}(z) - q_{k,-1}(z)
    + \sum_{n=1}^{2k-1}\frac{c_{n}}{z^{n}}.
\end{align*}
We let \( k=2 \),
use \eqref{eq:PPatalpha}
and then partial fraction decomposition
to find
that \( q \)
has the form
\begin{equation*}
q(z) = \frac{-1}{z-1} + \frac{1}{(z-1)^{2}}
	- \frac{1}{z-1} - \frac{1}{(z+1)^{2}}
	+ \sum_{n=1}^{3}\frac{c_{n}}{z^{n}}.
\end{equation*}
We substitute and find that the two relations
\eqref{eq:first_relation} and \eqref{eq:second_relation}
are satisfied if
\( c_{1}=2 \) and \( c_{2}=c_{3}=0 \).
Thus
\begin{equation*}
q(z) = \frac{-1}{z-1} + \frac{1}{(z-1)^{2}}
	- \frac{1}{z-1} - \frac{1}{(z+1)^{2}}
	+ \frac{2}{z},
\end{equation*}
is an RPF of weight \( 4 \) on \( G_{6} \)
with Hecke-symmetric ISP 
\( P_{\mathcal{A}_{2}} = \left\{ 1, -1 \right\} \).

\item 	
An RPF of weight \( 2k \)
with poles only in \( P_{\mathcal{A}_{3}}  \) 
has the form  
\begin{align*}
q(z)
& = q_{k,1/\sqrt{2}}(z) - q_{k,-\sqrt{2}}(z)
		+ \sum_{n=1}^{2k-1}\frac{c_{n}}{z^{n}}.
\end{align*}
For \( k=1 \)
we have that
\begin{equation}
\label{eq:RPF-G6hard-2}
q(z) = \frac{\sqrt{2}}{\sqrt{2}z-1} - \frac{1}{z+\sqrt{2}},
\end{equation}
is a rational period function of weight \( 2 \)
with the non-Hecke-symmetric ISP \( P_{\mathcal{A}_{3}} \).
\end{enumerate}
\end{enumerate}
\end{example}


\section{Families of rational period functions}

Schmidt and Sheingorn 
point out in \cite{MR1362251}
that conjugacy class generators
are \( \lambda_{p} \)-invariant.
By this they meant that
as \( p \) changes,
the values of the entries change
but as functions of \( \lambda \)
the generators do not change.
Rather,
increasing \( p \) 
changes the set of generators by adding matrices.
They called this phenomenon the ``\( q \)-principle.''
(Schmidt and Sheingorn used \( q \) to index the Hecke groups,
so in our context this would be the ``\( p \)-principle.'')

The \( p \)-principle
also holds for products of conjugacy class generators.
If a particular Hecke group \( G_{p_{0}} \)
has a conjugacy class with a product of generators
then 
that product 
(with each matrix a function of \( \lambda \))
is a conjugacy class block
for every Hecke group \( G_{p} \) with \( p \geq p_{0} \).
This gives us a natural method for constructing
families of RPFs across Hecke groups.
We illustrate this with an example.

\begin{example}
\label{ex:RPFfamily1}
Fix \( p \geq 3 \).
The product \( V_{1}V_{p-1} \) 
is a hyperbolic conjugacy class generator block
in \( G_{p} \).
Moreover,
\( (V_{1}V_{p-1})^{\top} = V_{1}V_{p-1} \)
so the corresponding ISP is Hecke-symmetric
and has \( 2 \) positive poles.
Now \[ V_{1}V_{p-1} = S(ST)^{p-2}S \sim S^{3}T(ST)^{p-3}, \]
so the corresponding reduced number
is 
\( [\overline{3,\underbrace{1, \ldots ,1}_{p-3}}] \)
and the set of positive poles in the ISP is
\begin{align*}
\mathcal{Z}_{\mathcal{A}} 
&  = \left\{ [2,\overline{\underbrace{1, \ldots ,1}_{p-3},3}], 
	[1,\overline{\underbrace{1, \ldots ,1}_{p-3},3}] \right\} 
 = \left\{ 
	\frac{\lambda_{p} + \sqrt{\lambda_{p}^{2}+4}}{2}, 
	\frac{-\lambda_{p} + \sqrt{\lambda_{p}^{2}+4}}{2} \right\}.
\end{align*}
This defines a family of ISPs
\( P_{\mathcal{A}} 
	= \mathcal{Z}_{\mathcal{A}} \cup T\mathcal{Z}_{\mathcal{A}} \)
for \( G_{p} \).
Then by \eqref{eq:RPFsymmetrickodd} we have that
\begin{align}
\label{eq:RPFfamily1}
q(z) 
& = \frac{1}{\left( z^{2} - \lambda_{p} z - 1 \right)^{k}}
  + \frac{1}{\left( z^{2} + \lambda_{p} z - 1 \right)^{k}}.
\end{align}
is an RPF of weight \( 2k \) (\( k \) odd) 
for any \( G_{p} \)
with poles 
\( P_{\mathcal{A}} = \mathcal{Z}_{\mathcal{A}} \cup T\mathcal{Z}_{\mathcal{A}} \).

If we allow \( p \) to take on any value
we have a family of ISPs 
\( P_{\mathcal{A}} = \mathcal{Z}_{\mathcal{A}} \cup T\mathcal{Z}_{\mathcal{A}} \)
and a corresponding family of RPFs for the Hecke groups.
This is the 
family 
constructed by Parson and Rosen \cite{MR748948}
and Schmidt \cite{MR1219337}.
Marvin Knopp's RPF \eqref{eq:MKsRPF}
for \( G_{3} \)
in Example \ref{ex:RPFsFewPoles} 
is a member of this family.

\end{example}


\subsection{Families of Hecke-symmetric ISPs}

The family of ISPs in Example \ref{ex:RPFfamily1}
is one of a larger class of Hecke-symmetric ISP families
for the conjugacy generator products
\( V_{j}V_{p-j} \).
For every fixed \( j \geq 1 \)
the product \( V_{j}V_{p-j} \)
produces a family of ISPs
for the Hecke groups with \( p \geq 2j+1 \).
The restriction on \( p \)
ensures that \( V_{j} \neq V_{p-j} \) 
and that the various families produced
are distinct.

We can produce other Hecke-symmetric classes of ISP families
by writing generator products that are conjugate to themselves.
We list blocks that produce other classes of families of ISPs.
\begin{enumerate}
\item 
Products of the form
\( V_{j}^{s}V_{p-j}^{s} \) for \( j, s \geq 1 \)
produce distinct families of ISPs
for the groups \( G_{p} \) with \( p \geq 2j+1 \).
\item 
Embedded products of the form in part 1, such as
\( V_{j}^{s}V_{\ell}^{t}V_{p-\ell}^{t}V_{p-j}^{s} \) 
for \( j, \ell, s, t \geq 1 \)
produce distinct families of ISPs
for the groups \( G_{p} \) with \( p \geq \max\{2j+1, 2\ell +1\} \).
\item 
The generator
\( V_{p/2} \) 
produces a family of ISPs for the groups \( G_{p} \) with \( p \) even.
\item 
The generator
\( V_{p/2} \) embedded in products of the form in parts \( 1 \) or \( 2 \), 
such as
\( V_{j}^{s}V_{p/2}V_{p-j}^{s} \) 
for \(j, s \geq 1 \)
produces distinct families of ISPs
for the groups \( G_{p} \) with \( p \) even, \( p \geq 2j+2 \).
\end{enumerate}

The next example gives two new families of RPFs 
obtained from two of these
self-conjugate generator products.
\begin{example}	
\label{ex:H-S-RPFfamilies}
\mbox{}
\begin{enumerate}
\item 	
Fix an even number \( p \) with \( p\geq 4 \).
Then \( V_{p/2} \) is a hyperbolic conjugacy class generator 
in \( G_{p} \)
that corresponds to a Hecke-symmetric ISP
with one positive pole.
Now 
\begin{equation*}
V_{p/2} = (ST)^{\frac{p}{2}-1}S
\sim S^{2}T(ST)^{\frac{p}{2}-2},
\end{equation*}
so the corresponding reduced number
is 
\[ \beta = [\overline{2,\underbrace{1, \ldots ,1}_{p/2-2}}], \]
and the positive pole is in
\begin{align*}
\mathcal{Z}_{\mathcal{A}} 
&  = \left\{ [1,\overline{\underbrace{1, \ldots ,1}_{p/2-2},2}] \right\} 
 = \left\{ 1 \right\} .
\end{align*}
Then \eqref{eq:RPFsymmetrickodd}
gives us that
\begin{align*}
q(z) 
& = \frac{1}{\left( z^{2} - 1 \right)^{k}}
\end{align*}
is a family of RPFs of weight \( 2k \) (\( k \) odd) 
for \( G_{p} \) (\( p \geq 4 \))
with poles 
\( P_{\mathcal{A}} = \mathcal{Z}_{\mathcal{A}} \cup T\mathcal{Z}_{\mathcal{A}} \).
The RPFs \eqref{eq:RPF-G4easy} for \( G_{4} \) 
and \eqref{eq:RPF-G6easy} for \( G_{6} \)
in Example \ref{ex:RPFsFewPoles} 
are members of this family.

\item 	
Fix \( p \geq 5 \).
Then the product \( V_{2}V_{p-2} \) 
is a primitive hyperbolic conjugacy class generator block
in \( G_{p} \).
The corresponding ISP is Hecke-symmetric
and has \( 2 \) positive poles.
Now \[ V_{2}V_{p-2} = STS(ST)^{p-3}S \sim S^{2}TS^{2}T(ST)^{p-4}, \]
so the corresponding reduced number
is 
\begin{equation*}
\beta = [\overline{2,2,\underbrace{1, \ldots ,1}_{p-4}}],\end{equation*}
and the set of positive poles in the ISP is
\begin{align*}
\mathcal{Z}_{\mathcal{A}} 
&  = \left\{ 
    [1,\overline{2,\underbrace{1, \ldots ,1}_{p-4},2}], 
	[1,\overline{\underbrace{1, \ldots ,1}_{p-4},2,2}] 
    \right\} \\
& = \left\{ 
	\frac{\lambda_{p}^{2}-2 + \sqrt{\lambda_{p}^{4}+4}}{2\lambda_{p}}, 
	\frac{-(\lambda_{p}^{2}-2) + \sqrt{\lambda_{p}^{4}+4}}{2\lambda_{p}}
    \right\}.
\end{align*}
If \( k \) is odd
\eqref{eq:RPFsymmetrickodd} gives us
\begin{align*}
q_{1}(z) 
& = \frac{1}{\left( \lambda_{p} z^{2}+(-\lambda_{p}^{2}+2)z-\lambda_{p} \right)^{k}}
  + \frac{1}{\left( \lambda_{p} z^{2}+(\lambda_{p}^{2}-2)z-\lambda_{p} \right)^{k}}.
\end{align*}
We multiply \( q_{1}(z) \) by \( \lambda_{p}^{k } \)
to conclude that
\begin{align}
\label{eq:RPFfamily3}
q(z) 
& = \frac{1}{\left( z^{2}+(-\lambda_{p}+2/\lambda_{p})z-1 \right)^{k}}
  + \frac{1}{\left( z^{2}+(\lambda_{p}-2/\lambda_{p})z-1 \right)^{k}},
\end{align}
is a family of RPFs of weight \( 2k \) (\( k \) odd) 
for \( G_{p} \) (\( p \geq 5 \))
with poles 
\( P_{\mathcal{A}} = \mathcal{Z}_{\mathcal{A}} \cup T\mathcal{Z}_{\mathcal{A}} \).

If we let \( p=3 \),
so \( \lambda_{p}=1 \) in \eqref{eq:RPFfamily3}
we have Marvin Knopp's RPF for the modular group \eqref{eq:MKsRPF},
even though the calculations above do not hold for \( p=3 \).
Thus we could think of \eqref{eq:MKsRPF}
as a member of both families
\eqref{eq:RPFfamily1} and \eqref{eq:RPFfamily3}.
\end{enumerate}
\end{example}


\subsection{Families of non-Hecke-symmetric ISPs}

We could construct families of non-Hecke-symmetric ISPs
\( P_{\mathcal{A}_{1}} \)
by starting with 
any conjugacy class generator product
that is not conjugate to its transpose.
We can use this
to write a family of RPFs with non-Hecke-symmetric poles,
at least for small weights.
We can also write a family of RPFs with poles 
in the Hecke-symmetric
\( P_{\mathcal{A}_{1}} \cup P_{\mathcal{A}_{2}} \),
where \( P_{\mathcal{A}_{2}} \) is the ISP 
conjugate to \( P_{\mathcal{A}_{1}} \).

\begin{example}		
\label{ex:non-H-S-RPFfamilies}
Fix \( p \geq 5 \).
Then \( V_{2} \) 
is a hyperbolic conjugacy class generator 
in \( G_{p} \)
that is not conjugate to its transpose.
(\( V_{2} \) is parabolic in \( G_{3} \)
and equal to its transpose in \( G_{4} \).)
The corresponding ISP is non-Hecke-symmetric
and has one positive pole.
Now \( V_{2} = STS \sim S^{2}T \)
so the corresponding reduced number
is 
\( [\overline{2}] \),
and 
\( \mathcal{Z}_{\mathcal{A}_{1}} 
  = \left\{ [1,\overline{2}] \right\} 
  = \left\{ \sqrt{\lambda_{p}^{2}-1} \right\} \).
The ISP is 
\begin{equation*}
P_{\mathcal{A}_{1}} 
	= \left\{ \sqrt{\lambda_{p}^{2}-1}, {-1}/{\sqrt{\lambda_{p}^{2}-1}} \right\}.
\end{equation*}

The conjugacy class generator
\( V_{2}^{\top} = V_{p-2} \)
gives us the ISP conjugate to \( P_{\mathcal{A}_{1}} \).
Now \( V_{p-2} = (ST)^{p-3}S \sim S^{2}T(ST)^{p-4}, \)
so the corresponding reduced number
is 
\( [\overline{2;\underbrace{1, \ldots ,1}_{p-4}}] \),
and
\( \mathcal{Z}_{\mathcal{A}_{2}} 
  = \left\{ [1;\overline{\underbrace{1, \ldots ,1}_{p-4},2}] \right\} 
  = \left\{ {1}/{\sqrt{\lambda_{p}^{2}-1}} \right\} \).
The ISP is 
\begin{equation*}
P_{\mathcal{A}_{2}} 
	= \left\{ {1}/{\sqrt{\lambda_{p}^{2}-1}}, -\sqrt{\lambda_{p}^{2}-1} \right\},
\end{equation*}

By \eqref{eq:RPFnonsymmetrick}
an RPF of weight \( 2k \)
with poles only in \( P_{\mathcal{A}_{1}}  \) 
has the form  
\begin{align*}
q(z)
& = q_{k,\sqrt{\lambda_{p}^{2}-1}}(z) - q_{k,-1/\sqrt{\lambda_{p}^{2}-1}}(z)
	+ \sum_{n=1}^{2k-1}\frac{c_{n}}{z^{n}}.
\end{align*}
We let \( k=1 \)
and substitute to find that the two relations
\eqref{eq:first_relation} and \eqref{eq:second_relation}
are satisfied if \( c_{1} = 0 \).
This gives us that
\begin{equation*}
\label{eq:RPF-family-hard-1}
q(z) = \frac{1}{z-\sqrt{\lambda_{p}^{2}-1}} 
	- \frac{\sqrt{\lambda_{p}^{2}-1}}{\sqrt{\lambda_{p}^{2}-1}z+1},
\end{equation*}
is a family of RPFs of weight \( 2 \) 
for \( G_{p} \) (\( p \geq 5 \))
with poles 
in the non-Hecke-symmetric ISP
\( P_{\mathcal{A}_{1}} \).
The RPFs \eqref{eq:RPF-G5hard-1} for \( G_{5} \) 
and \eqref{eq:RPF-G6hard-1} for \( G_{6} \)
in Example \ref{ex:RPFsFewPolesHard} 
are members of this family.

In a similar fashion
we can use \eqref{eq:RPFnonsymmetrick}
to find that
\begin{equation*}
q(z) = \frac{\sqrt{\lambda_{p}^{2}-1}}{\sqrt{\lambda_{p}^{2}-1}z-1} 
	- \frac{1}{z+\sqrt{\lambda_{p}^{2}-1}},
\end{equation*}
is a family of RPFs of weight \( 2 \) 
for \( G_{p} \) (\( p \geq 5 \))
with poles 
in the non-Hecke-symmetric ISP
\( P_{\mathcal{A}_{2}} \).
The RPFs \eqref{eq:RPF-G5hard-2} for \( G_{5} \) 
and \eqref{eq:RPF-G6hard-2} for \( G_{6} \)
in Example \ref{ex:RPFsFewPolesHard} 
are members of this family.

Finally, we can write a family of RPFs
with poles in the Hecke-symmetric
\( P_{\mathcal{A}_{1}} \cup P_{\mathcal{A}_{2}} \).
We use 
\eqref{eq:RPFquadraticanyk}
to get that
\begin{equation*}
q(z) = \frac{1}{(z^{2}-(\lambda_{p}^{2}-1))^{k}} 
	-\frac{(-1)^{k}}{((\lambda_{p}^{2}-1))z^{2}-1)^{k}},
\end{equation*}
is a family of RPFs of weight \( 2k \) (for any \( k \))
in \( G_{p} \) (\( p \geq 5 \))
with poles in
\( P_{\mathcal{A}_{1}} \cup P_{\mathcal{A}_{2}} \).
The RPFs \eqref{ex:RPF-G5easy} for \( G_{5} \) 
and \eqref{ex:RPF-G6easy-nonHS} for \( G_{6} \)
in Example \ref{ex:RPFsFewPolesHard} 
are members of this family.
\end{example}

\bibliographystyle{alpha}

\bibliography{WRLibraryMR,WRLibraryNotMR}

\end{document}